\documentclass[12pt]{article}

\usepackage{amsmath}
\usepackage{amssymb}
\usepackage{amsfonts}
\usepackage{epsfig}

\newtheorem{theo}{Theorem}

\newtheorem{prop}{Proposition}
\newtheorem{lemm}{Lemma}

\newcommand{\eqco}{\setcounter{equation}{0}}
\newcommand{\thco}{\setcounter{theo}{0}}
\newcommand{\prco}{\setcounter{prop}{0}}
\newcommand{\laco}{\setcounter{lemm}{0}}
\newcommand{\coco}{\setcounter{coro}{0}}
\newcommand{\deco}{\setcounter{defn}{0}}
\newcommand{\allco}{\eqco  \thco \prco \laco \coco \deco}
\def\N{\mathbb{N}}
\def\Z{\mathbb{Z}}
\def\R{\mathbb{R}}
\def\E{\mathbb{E}\,}

\newcommand{\la}{\lambda}
\newcommand{\ka}{\kappa}
\newcommand{\K}{{\cal K}}
\newcommand{\tolone}{\stackrel{{L^1}}{\longrightarrow}}

\newcommand{\tod}{\stackrel{{\cal D}}{\longrightarrow}}
\newcommand{\toas}{\stackrel{{\rm a.s.}}{\longrightarrow}}
\newcommand{\eqd}{\stackrel{{\cal D}}{=}}
\newcommand{\toP}{\stackrel{{P}}{\longrightarrow}}
\newcommand{\toone}{\stackrel{{L^1}}{\longrightarrow}}
\newcommand{\toLp}{\stackrel{{L^p}}{\longrightarrow}}
\newcommand{\limn}{\lim_{n \to \infty}}
\newcommand{\Po}{{\cal P}}
\newcommand{\Tau}{{\cal T}}
\newcommand{\tPo}{\tilde{{\cal P}}}
\newcommand{\tdelta}{\tilde{\Delta}}
\newcommand{\B}{{\cal B}}
\newcommand{\Ed}{{\cal E}}
\newcommand{\A}{{\cal A}}
\newcommand{\Op}{{\cal O}}
\newcommand{\RR}{{\cal R}}
\newcommand{\MST}{{\rm MST}}

\newcommand{\W}{{\cal W}}
\newcommand{\diam}{{\rm diam}}
\newcommand{\G}{{\cal G}}
\newcommand{\HH}{{h}}
\newcommand{\HHH}{{\zeta}}
\newcommand{\bb}{{\bf b}}
\newcommand{\balpha}{{\bf \alpha}}
\newcommand{\0}{{\bf 0}}
\newcommand{\X}{{\cal X}}
\newcommand{\YY}{{\cal Y}}
\newcommand{\TT}{{T}}
\newcommand{\U}{{\cal U}}

\newcommand{\Y}{{\bf Y}}

\newcommand{\cc}{{\rm c}}

\newcommand{\NN}{{\cal N}}
\newcommand{\F}{{\cal F}}

\newcommand{\tlB}{\widetilde{tB}_0}
\newcommand{\tnlB}{B_n}

\newcommand{\eps}{\varepsilon}
\newcommand{\bdm}{\begin{displaymath}}
\newcommand{\edm}{\end{displaymath}}
\newcommand{\be}{\begin{equation}}
\newcommand{\ee}{\end{equation}}

\newcommand{\bea}{\begin{eqnarray}}
\newcommand{\eea}{\end{eqnarray}}
\newcommand{\bear}{\begin{eqnarray}}
\newcommand{\eear}{\end{eqnarray}}
\newcommand{\bean}{\begin{eqnarray*}}
\newcommand{\eean}{\end{eqnarray*}}

\newcommand{\card}{{\rm card}}

\newcommand{\Cov}{{\rm Cov}}
\newcommand{\cov}{{\rm Cov}}
\newcommand{\Var}{{\rm Var}}

\newcommand{\qed}{\rule[-1mm]{3mm}{3mm}}

\renewcommand{\preceq}{\preccurlyeq}

\newcommand{\delinf}{\delta_\infty(\Po_\la)}
\newcommand{\delinfi}{\delta_\infty^i(\Po_\la)}
\newcommand{\delinfj}{\delta_\infty^j(\Po_\la)}
\newcommand{\delinfd}{\delta'_\infty(\Po_\la)}

\begin{document}
\title{\bf 
Multivariate spatial central limit theorems
with applications to percolation and spatial graphs
}

\author{MATHEW D.~PENROSE
\\
University of Bath  \\
} 
%
\maketitle

\begin{abstract}
Suppose $X = (X_x, x$ in $Z^d)$ is a family of i.i.d. 
variables in some measurable space,
$B_0$ is a bounded set in $R^d$, 
and for $t > 1$, $H_t$ is a measure on $tB_0$ 
  determined by the restriction of $X$ to  
lattice sites in or adjacent to
$tB_0$. 
We prove convergence to a white noise process for  
the random measure on $B_0$ 
given by
$t^{-d/2}(H_t(tA)-EH_t(tA))$ for subsets $A$ of $B_0$, as $t$ becomes large,
 subject to $H$ satisfying a ``stabilization''
condition (whereby the effect of changing $X$ at a single site $x$
is local) but with no assumptions on the rate of decay of
correlations.  We also give a multivariate central limit theorem
for the joint distributions of two or more such measures
$H_t$, and adapt the result to measures based
on Poisson and  binomial point processes.
Applications given include a white noise limit for the measure which
counts  clusters  of critical percolation, a functional central limit
theorem for the empirical process of the edge lengths of the minimal
spanning tree on random points, and central limit theorems
for the on-line nearest neighbour graph.
\end{abstract}

Short title: {\bf Multivariate spatial CLTs} \\

{\em Key words and phrases.}
Central limit theorem, white noise, minimal spanning tree,
empirical process, on-line nearest neighbour graph, percolation. \\

{\em American Mathematical Society 2000 Classifications.} 
Primary-60F05, 60D05;
Secondary-05C80, 60K35. 

\newpage


\section{Introduction}
Several approaches have been developed  for
proving central limit theorems for random variables
which arise as the sum of contributions from points of a Poisson
or binomial point process in $\R^d$, when each contribution
is locally determined in some sense. These include
Stein's method, the method of moments,
and a martingale method. 
%
%

By keeping track of the location of each contribution in $\R^d$,
one can often in a natural way create a random measure, i.e. a random field
indexed by subsets of $\R^d$ or by test functions on $\R^d$.
 It  is of interest to look for 
 multivariate central limit theorems for such random fields,
typically with weak convergence of finite-dimensional distributions
to those of white noise. 
Multivariate central 
limit theorems of this type were recently derived using the method
of moments by Baryshnikov and Yukich \cite{BY} and can also be derived using 
Stein's method \cite{PY6}. 
Both of these methods seem to require, in addition to a `stabilization'
condition which formalizes the locally determined contributions,
a form of exponential decay of spatial correlations.


It is of interest to extend  these results
to cases which satisfy stabilization but are not believed to
satisfy  exponential decay.
These include, for example, measures  associated 
with the  minimal spanning tree (MST) on a Poisson point process
or with critical percolation.
The martingale method is especially powerful 
in giving central limit theorems for these examples (Kesten and Lee
\cite{KL}, Zhang \cite{Z}).
For an exposition of this method in a general setting, see Penrose \cite{Pclt},
Penrose and Yukich \cite{PY}.
However, these works do not address the convergence to
white noise of random measures.

In the present paper we extend the martingale method to
give such a convergence to white noise of  stabilizing
 random fields indexed by subsets of $\R^d$,
 and illustrate the method both for percolation
 and for the minimal spanning tree.   
A further multivariate  direction in which we extend the existing literature
is  by considering convergence to a multivariate normal for
two or more random fields based on the same underlying
spatial process; in particular,
we shall show that the finite-dimensional distributions of the
empirical process
of the lengths of the MST on random points, suitably scaled
and centred, converge to those of a certain Gaussian process.

A further aim of this paper is to  treat discrete examples
(such as percolation) and continuous ones (such as
the MST) in a unified manner. In this spirit
we shall derive our basic general result for Poissonian continuum 
systems (Theorem \ref{1Hclt})  by direct application
of the basic result for lattice systems (Theorem \ref{basictheo}),
although some extra work is needed to give the limiting
covariances for Theorem \ref{1Hclt} in a more explicit
form than was done in previous continuum central limit
theorems proved by the martingale method.
For this reason, in the continuum we consider only uniform
densities of points over a fixed $d$-dimensional set
(denoted $B_0$ in the sequel), unlike Baryshnikov and Yukich
who consider non-uniform densities. It is likely that
with extra work, the martingale-based proof of 
 Theorem \ref{1Hclt} could be extended to give
 multivariate  central
limit theorems  for  point processes with non-uniform
densities.  For martingale-based proofs of univariate
central limit theorems on non-uniform points,
see Lee \cite{LeeII} and Penrose (\cite{Pbk}, Section 13.7).

In the continuum, stabilizing random fields are often defined in terms of
   graphs $G$ which are themselves stabilizing, i.e., locally determined
in a certain sense.  Stabilizing graphs include
 the  MST, $k$-nearest neighbour, and sphere of influence graphs.
 Given  a stabilizing graph $G$,
the theory presented here applies to random fields
(indexed by subsets $A$ of $\R^d$)
which count, for example, 
 the number of leaves of $G$ in $A$, the number
 of components of $G$ that include vertices in $A$, or the
 sum of weighted edge lengths
$\phi(\vert e \vert)$ with the sum over edges $e$ of
$G$ having endpoints in $A$.  

One stabilizing graph which has not been considered in previous
discussions of stabilizing graphs is the {\em on-line nearest
neighbour graph}, in which  random points in $B_0$ are 
 randomly ordered, and each point (except the first) is
 connected to its nearest neighbour amongst its predecessors
in the ordering. This graph is of recent interest in
connection with the modelling of scale-free networks
\cite{BBBCR,FKP}. Unlike methods based on exponential
decay, our methods provide
 central limit theorems for this graph too; see Section \ref{secOL}.

The rest of the paper is laid out as follows. The next section
contains statements of the main general results. 
Section \ref{Secaps} contains applications of these to
percolation, MST, and nearest-neighbour type graphs.
 Sections \ref{secprfdscrt}
 and \ref{secprfcont}
contain proofs of the general results.

\section{General multivariate central limit theorems}
\allco

\subsection{Notation used throughout}
\label{secnotation}

Let $d \geq 1$ be an integer and let $\0$ denote the origin of
$\R^d$. For $x \in \R^d$, write $|x|$ for the Euclidean norm
of $x$.
For $A \subseteq \R^d$, $t \in \R$, 
 and $y \in \R^d$, 
let $tA$ denote the scaled  set $\{tx: x \in A\}$, and
let $\tau_y(A)$ denote the translated set $\{y+x: x \in A\}$.
Let $\partial(A)$ denote the boundary of $A$, that is, 
the intersection of the closure of $A$ with that of
its complement.
If $A$ is (Lebesgue) measurable, write $|A|$ for
its Lebesgue measure, 
and if $A$ is finite, write $\card(A)$ for
the number of elements of $A$.
Write $\diam(A)$ for $\sup\{|x-y|:x\in A, y \in A\}$.
Given a sequence of sets $(A_n)_{n \geq 1}$,
 write $\liminf(A_n)$ for $\cup_{n=1}^{\infty}
\left( \cap_{m=n}^\infty
A_m \right)$. 

For $z \in \Z^d$, and $\eps >0$, 
let $Q_z^\eps$ denote the cube $\tau_{\eps z}([-\eps,0)^d)$.
For $x \in \R^d$, and $r>0$, let $B_r(x)$ be the closed Euclidean
ball of radius $r$ centred at $x$. Let $\rho$
be a finite constant, satisfying
$\rho  \geq \sqrt{d}$ but otherwise arbitrary.
 For $A \subseteq \R^d$, 
let $\widetilde{A}$ denote the discretization of $A$ given by
\bea
\widetilde{A}:= \{z \in \Z^d: B_\rho(z) \cap A \neq  \emptyset \}.
\label{1121}
\eea
The condition $\rho \geq \sqrt{d}$ ensures that
 $Q_z^\eps \subset B_\rho(\eps z)$
for $\eps \leq 1$.

For $\sigma > 0$, let 
 $\NN(0,\sigma^2)$  be the normal probability distribution  on $\R$ with
density 
$f(x) =
 (2 \pi \sigma^2)^{-1/2}
 \exp( - x^2/(2 \sigma^2))$.
Also,  let 
$\NN(0,0)$ represent the degenerate probability distribution on
$\R$ consisting of 
a unit point mass at zero, which we view as a special case of the normal.
Given a nonnegative definite $k \times k$ matrix $\Sigma$, 
let $\NN(\0, \Sigma)$ denote
the centred multivariate normal distribution with
covariance matrix $\Sigma$, i.e. the distribution of
a random $k$-vector ${\bf X}$ satisfying 
${\bf a'} {\bf X}
 \sim \NN(0,{\bf a'} \Sigma {\bf a})$ for
all deterministic $k$-vectors ${\bf a}$
(this definition includes the case
when $\Sigma$ is singular).
Denote convergence in probability by $\toP$,
convergence in $p$th moment by $\toLp$,  convergence in distribution
by $\tod$,
and  denote equality in distribution by $\eqd$.

We say a subset of $\R^d$ is  {\em Riemann measurable}
if it has Riemann integrable indicator function; in other words
(see Rudin \cite{Rud}),  
we say a subset of $\R^d$ is   Riemann measurable
 if it is bounded and has Lebesgue-null boundary. 
Let $\RR(\R^d)$ denote the collection of Riemann measurable subsets of $\R^d$.
In the sequel, we shall assume $B_0$ is a fixed set in $\RR(\R^d)$ (so
in particular $B_0$ is bounded); we shall also assume that
$|B_0| >0$, which is equivalent to assuming that
$B_0$ has non-empty interior.
For example,
$B_0$ could be the $d$-dimensional unit cube.
Let $\RR(B_0)$ denote the collection of Riemann measurable
 subsets of $B_0$.

For $y,z \in \Z^d$, write $y \prec z$ if $y$ precedes $z$
in the lexicographic ordering on $\Z^d$, and $y \preceq z$
if either $y \prec z$ or $y= z$.

\subsection{A central limit theorem for lattice systems}
\label{secmainth}

Let $(E,{\cal E},P_0)$ be an arbitrary probability space.
On a suitable probability space $(\Omega, \F, P)$, let
 $X = (X_z, z \in \Z^d)$ be 
be a family of independent identically distributed
random elements of $E$,
each $X_z$ having distribution $P_0$,
  indexed by the integer lattice, and let $X_*$ be 
a further
$E$-valued variable with distribution $P_0$, independent of $X$
(i.e.,  an independent copy
of $X_\0$).
For existence of such an $(\Omega,\F,P) $ and $X$, see for example
section 8.7 of Williams \cite{W}.
For $y \in \Z^d$, let 
$\tau_yX$ denote the translated family
of variables $(X_{z+y},z \in \Z^d)$.

Suppose $B_0 \in \RR(\R^d)$ with $|B_0|>0$.
By a {\em random set function} on $B_0$ we mean a collection
$H= (H_t(A): t \geq 1, A \in \RR(B_0))$, where for
each $t \geq 1$ and $A \in \RR(B_0)$,  $H_t(A)$ is
a random variable
that is a function of $(X_z, z \in \tlB)$, so that,
strictly speaking,  $H_t(A)$ is itself 
a measurable function from $E^{\tlB}$ to $\R$.
If we wish to emphasize the dependence on $X$  of the
value of $H_t(A)$
we write
$H_t(X,A)$ for $H_t(A)$. In many examples
 $H_t(\cdot)$ is a (random)  measure or outer measure
(see, e.g., Durrett \cite{Dur2})
on Borel subsets of
$tB_0$ but we do not need to assume this for the general result
(we restrict attention to Riemann measurable $A$).
For $t \geq 1$, $y \in \Z^d$, and  $A \in \RR(B_0)$, define
\bea
H_{t,y}(A) = H_{t,y}(X,A) := H_t (\tau_{y} X,A)
\label{Htzdef}
\eea
In \cite{Pclt}, a CLT   is established under a 
  stabilization condition
 which says, loosely speaking, that 
 the effect on a given random set function $H$ of resampling 
 the value of $X$ at a single site is local. 
To extend these to central limit theorems for random fields,
we require a 
 modification of the
stabilization condition used in \cite{Pclt}.

Let $X^\0$ be the process $X$ with the value $X_\0$ at the origin replaced by 
the independent copy $X_*$ of $X_\0$, but with the values at
 all other sites the same (i.e., $X^\0 = (X^\0_z, z \in \Z^d)$ 
is the random field given by $X^\0_\0 = X_*$ and
$X^\0_z = X_z$ for $z \in \Z^d \setminus \{\0\}$). 
Given a random set function $H$,
 define the increment $\Delta_{t,y}^H(A)$
for $y \in \Z^d$ and $ t \geq 1, A \in \RR(B_0)$,
by
\bea
\Delta_{t,y}^H(A) = H_{t,y}(X,A) - H_{t,y}(X^\0,A).
\label{deltadef}
\eea
We consider random set functions $H$ with the property that
there exists a random variable $ \Delta^H_\infty$ such that for 
 all $A \in \RR(B_0)$, and all
$[1,\infty) \times \Z^d$-valued sequences 
  $(t_n,y_n)_{n \geq 1}$:
\bea
\Delta_{t_n,y_n}^H (A) \toP \Delta^H_\infty ~~~~{\rm 
if} ~~ \liminf_{ n \to \infty} (\tau_{y_n}(t_n A ) )=  \R^d
\label{0414a}
\eea
and also
\bea
\Delta_{t_n,y_n}^H (A) \toP 0 ~~~~{\rm 
if} ~~ \liminf_{ n \to \infty} (\tau_{y_n}(t_n (B_0 \setminus A) ) )=  \R^d.
\label{0414b}
\eea 
Eqns (\ref{0414a}) and (\ref{0414b}) are our stabilization conditions. 
The second condition (\ref{0414b}) is a novel feature of this paper;
it was not required for 
the  CLTs presented in \cite{Pclt,PY}. 
The first condition (\ref{0414a})
is similar to the stabilization
condition in Definition 2.3 of \cite{Pclt}.

We shall require also that there exist $\gamma >2$ such that 
the moments condition
\begin{equation}
\sup\{ \E[|\Delta_{t,-y}^H(A) |^\gamma] :  
A \in \RR(B_0), t \geq 1, y \in \tlB \}  < \infty
\label{4moments}
\end{equation}
is satisfied. Observe that 
$\Delta_{t,-y}^H(A)  $ is identically zero  for $y \in
\Z^d \setminus \tlB$, and therefore 
condition (\ref{4moments}) is equivalent to
\bea
\sup\{ \E[|\Delta_{t,-y}^H(A) |^\gamma] :  
A \in \RR(B_0), t \geq 1, y \in \Z^d \}  < \infty.
\label{4moments2}
\eea

For $y \in \Z^d$, let $\F_y$ be the $\sigma$-field generated by
 $(X_z,z\preceq  y)$ (recall that $\preceq$ denotes
the lexicographic ordering on $\Z^d$).
Now we can state our main general result for lattice systems.

\begin{theo}
\label{basictheo}
Suppose $B_0 \in \RR(\R^d)$ with $|B_0|>0$.
Suppose that $H^1, \ldots, H^k$ are random set functions on $B_0$,
each of which satisfies the 
 stabilization conditions $(\ref{0414a})$ and $(\ref{0414b})$,
along with the moments condition $(\ref{4moments})$ for
some $\gamma >2$.
Let the $k \times k$ matrix $(\sigma^*_{ij})_{i,j=1}^k $ be given by
\bea
\sigma^*_{ij} := 
\E[ \E(\Delta^{H^i}_\infty|\F_\0) \E(\Delta^{H^j}_\infty|\F_\0)  ] .
\label{0905a}
\eea
Then if $A_,\ldots,A_k$ are
 Riemann measurable subsets of $B_0$,
 for $1 \leq i \leq j \leq k$ 
we have
\begin{equation}
\lim_{t \to \infty} t^{-d} \Cov(H_t^i(A_i), H_t^j(A_j)) = 
\sigma^*_{ij} |A_i \cap A_j|,
\label{limvar}
\end{equation}
and as $t \to \infty$,
\begin{equation}
(t^{-d/2} ( H_t^i(A_i) - \E H_t^i(A_i) ))_{i=1}^k \tod 
\NN(\0, (\sigma^*_{ij} |A_i \cap A_j|)_{i,j=1}^k ).
\label{limdist}
\end{equation}
\end{theo}
In many examples, we consider only the case of a single
random set function, that is, the case where
each of $H^1,\ldots, H^k$ are all the same random set function $H$. In this
case the result says  that all the finite-dimensional  joint distributions 
of 
 $ (t^{-d/2} ( H_t(A) - \E H_t(A)), A \in \RR(B_0) )$,
converge to  those of a centred Gaussian
 process $(W(A), A \in \RR(B_0))$ with covariance function 
$$
\E [W(A)W(A') ] = |A \cap A'|
\E[( \E[\Delta_\infty^H|\F_\0])^2 ],
$$
 i.e., a white noise process.

\subsection{Central limit theorems for continuum systems}
\label{contsubsec}

By a {\em point process set function} we mean
a real-valued functional $\HH(\X,A)$  
defined for all $A\in \RR(\R^d)$ and
  finite subsets $\X$ of $\R^d$, such that
\begin{enumerate}
\item[(i)]
$(x_1,\ldots,x_k) \mapsto h(\{x_1,\ldots,x_k\},A)$
is a Borel-measurable function, for all 
$k \in\N, A \in \RR(\R^d)$; 
\item[(ii)]
for all $A \in \RR(\R^d),$  $ y \in \R^d,$
and all finite  $ \X \subset \R^d$,
$h$ satisfies the
translation-invariance condition
\bea
\HH(\tau_y(\X ),\tau_y( A )) = \HH(\X,A).
\label{transinv}
\eea
\end{enumerate}

For $\lambda >0$, let $\Po_\la$ denote a homogeneous
Poisson point process in $\R^d$ of intensity $\lambda$
(viewed as a random subset of $\R^d$). 
Given $B_0 \in \RR(\R^d)$ with $|B_0| >0$,
 define the point processes
$$
\Po_{\la,t} :=  \Po_\la \cap (tB_0),~~~~ t \geq 1.
$$
 We derive a multivariate  central limit theorem for
$(\HH(\Po_{\la,t},tA), A \in \RR(B_0))$ as $t \to \infty$.
The conditions on $\HH$  for our central limit theorem are defined
in terms of  the ``add one cost on $A$'' defined by
\bea
\delta(A,\X) := \HH(\X \cup \{\0\},A) - \HH(\X,A).
\label{add1c}
\eea
We shall say the point process set function $h$ 
is {\em strongly stabilizing} at intensity $\la$ 
if there exist almost surely finite  random variables
$\delinf$ 
(the {\em stabilizing limit} of $h$ at intensity $\la$)
and  $S$ (a {\em radius of
stabilization} of $\HH$ at intensity $\la$) such that with probability 1,
$\Po_\la$ is such that for $A \in \RR(\R^d)$ and for all finite
${\cal A} \subset (\R^d \setminus B_S(\0))$, 
\bea
 \delta( A ,(\Po_\la \cap {B_S(\0)}) \cup {\cal A} ) = \delinf
~~~
{\rm if} ~~A \supseteq B_S(\0),
\label{sstab1}
\eea
and
\bea
 \delta( A ,(\Po_\la \cap {B_S(\0)}) \cup {\cal A} ) = 0 
~~~ {\rm if} ~~A \cap B_S(\0) = \emptyset.
\label{sstab2}
\eea
Thus, $S$ is a radius of stabilization if the add one cost on $A$ for
 the restriction of $\Po_\la$
to a region containing the ball $B_S(\0)$
 is unaffected by changes in the configuration
outside the ball $B_S(\0)$ if 
$B_S(\0) \subseteq A$
or   $B_S(\0) \cap A = \emptyset$,
taking the value $\delinf$  if
$B_S(\0) \subseteq A$
and the value zero if $B_S(\0) \cap A = \emptyset$.
Our notion of strong stabilization (\ref{sstab1}) is similar to
 that used in \cite{PY}.
The second stabilization condition (\ref{sstab2}),
like its discrete counterpart (\ref{0414b}), is new to this paper.  

As in \cite{PY}, as well as strong stabilization 
we have a notion of `weak stabilization' which we
shall describe in Section \ref{secprfcont}. Loosely speaking,
the distinction is that in (\ref{sstab1}) and (\ref{sstab2}),
the set ${\cal A}$ runs through {\em all} finite
sets in $\R^d \setminus B_S(\0)$, whereas the
corresponding weak stabilization
conditions (eqns (\ref{0513a}) and (\ref{0513b})  
below) refer only to subsets of the underlying
Poisson process $\Po_\la$. Theorem \ref{1Hclt} 
below is stated under the strong stabilization conditions
(\ref{sstab1}) and (\ref{sstab2}) but actually still holds
if these are replaced by the weak stabilization
conditions (\ref{0513a}) and (\ref{0513b}).
Theorems \ref{2Hclt} and \ref{Xncoro} really require
the strong stabilization conditions (\ref{sstab1}) and (\ref{sstab2}).
All the examples discussed here satisfy (\ref{sstab1}) and (\ref{sstab2})
but there may be examples satisfying weak but not
strong stabilization, for example in relation to
germ-grain (Boolean) models \cite{BY,HM,MRbk} with 
no bound on grain sizes, or to
the random connection model with long-range connections
\cite{MRbk,RS}.

Let $\la >0$ and let 
$B_0 \in \RR(\R^d)$ with $|B_0|>0$.
Given $t \in [1,\infty)$ and $m \in \N$,
 let $\U_{m,t}$ be a point process
consisting of $m$ independent random $d$-vectors, each of
them uniformly distributed on $tB_0$.
Also let $\mu_{\la,t} $ be the expected number of points of $\Po_{\la,t}$, i.e., let
\bea
\mu_{\la,t} : = \lambda t^d |B_0|
\label{mudef}
\eea
 We consider  functionals $\HH$ satisfying the  moments
condition
\bea
\sup_{t \geq 1, A \in \RR(B_0),
 x \in -tB_0}
~~~
 \sup_{m \in [\mu_{\la,t}/2,3\mu_{\la,t}/2]} \{
\E [   \delta(\tau_x(tA),\tau_x({\cal U}_{m,t}) ) ^4 ] \} < \infty.
\label{ubdmom}
\eea
In the sequel, it is likely that the fourth moments condition (\ref{ubdmom})
  can be replaced by a $2 +
\epsilon$ moment condition, but this would not greatly expand
the range of applications known to the author.

We also require  a mild uniform
 bound on $\HH$ in terms of the size of  $\X$, whereby there
exists a constant $\beta_2$ such that for all finite sets $\X
\subset \R^d$, and all $A\in \RR(\R^d)$,
\bea
| \HH(\X,A ) | \leq \beta_2 (\diam (\X) + \card (\X))^{\beta_2}.
\label{polybd}
\eea

We now give a multivariate CLT for  $h(\Po_{\la,t},tA )$.
This is the first of our continuum analogues to Theorem \ref{basictheo}.
\begin{theo}
\label{1Hclt}
Let $\la >0$, and let
$B_0 \in B(\R^d)$ with $|B_0|>0$.
  Suppose that $\HH^1, \ldots, \HH^k$ 
are point process set functions 
 which satisfy the stabilization conditions
 $(\ref{sstab1})$, $(\ref{sstab2})$, 
 the moments condition $(\ref{ubdmom})$, and the uniform bound
 $(\ref{polybd})$.
 Define the $k \times k$ matrix $(\sigma_{ij}^\la)_{i,j=1}^k $  by
\bea
\sigma_{ij}^\la=  \E[ \E(\delinfi| \F) \E(\delinfj|\F)],
\label{0905}
\eea
where $\F$ denotes the $\sigma$-field generated by
the Poisson configuration in the half-space $\{x = (x_1, \ldots, x^d) \in
\R^d: x_1 < 0\}$, and $\delinfi$ is the stabilizing limit of
$h^i$. 
If $A_1,\ldots, A_k$ are sets in $\RR(B_0)$, then
as $t \to \infty$,  
\bea
t^{-d} \cov(h^i(\Po_{\la,t},tA_i), h^j(\Po_{\la,t},tA_j))
\to \lambda \sigma_{ij}^\la|A_i \cap A_j|
\label{0516a}
\eea
 and
\bea
t^{-d/2} (\HH^i(\Po_{\la,t},A_i) - \E \HH^i( \Po_{\la,t},A_i) )_{i=1}^k \tod 
\NN(\0,  (\lambda \sigma_{ij}^\la|A_i \cap A_j|)_{i,j=1}^k).
\label{0516b}
\eea
\end{theo}


The next result is a de-Poissonized version of  Theorem \ref{1Hclt}.

\begin{theo}
\label{2Hclt}
Let $\la >0$ and suppose 
$B_0 \in \RR(\R^d) $ with $|B_0|>0$.
   Suppose that $\HH^1, \ldots, \HH^k$ are point process set functions,
each $h^j$ satisfying the strong stabilization conditions
$(\ref{sstab1})$ $($with stabilizing limit denoted $\delinfj$)
 and $(\ref{sstab2})$, along with the
 moments condition $(\ref{ubdmom})$, and
the uniform bound $(\ref{polybd})$.
Let the matrix
 $(\sigma_{ij}^\la )_{i,j=1}^k$ 
be given by
$(\ref{0905})$.
Then if $A_1, \ldots, A_k$ are sets in $\RR(B_0)$,   
if we define the matrix 
$\Tau^\la := (\tau_{ij}^\la)_{i,j=1}^k$ by
\bea
\tau_{ij}^\la  :=  \frac{\sigma_{ij}^\la |A_i \cap A_j|}{|B_0|} 
- \frac{|A_i|\cdot|A_j|}{|B_0|^2}
\E [\delinfi] \E [\delinfj],
\label{taudef}
\eea
and if  $(t_n)_{n \geq 1}$ is a $[1,\infty)$-valued sequence satisfying 
\bea
\label{tncond}
\limsup_{n \to \infty} n^{-1/2} |(\la t_n^d |B_0| -n) | < \infty, 
\eea
then for $i, j \in \{1,\ldots, k\}$, 
\begin{equation}
\limn n^{-1} \cov(\HH^i(\U_{n,t_n},t_n A_i), \HH^j(\U_{n,t_n},t_n A_j ))
 =   \tau_{ij}^\la
\label{varlim}
\end{equation}
 and
 as $n \to \infty$,
\begin{equation}
n^{-1/2} (\HH^i(\U_{n,t_n},t_nA_i) - \E \HH^i(\U_{n,t_n},t_nA_i) )_{i=1}^k
 \tod \NN(\0, \Tau^\la).
\label{2clteq}
\end{equation}
\end{theo}

 Given $\gamma \in \R$, we  shall say $\HH$
is {\em homogeneous of order $\gamma$} if
\bea
\HH(a \X,aA) = a^\gamma \HH(\X,A), ~~~~\forall a \in \R,~A \in \RR(\R^d),
{\rm ~finite}~\X \subset \R^d.
\label{eqhomog}
\eea
If $\HH$ satisfies homogeneity, it is easy to deduce from the above
theorems a multivariate CLT, either for a homogeneous Poisson processes of 
intensity $\lambda$  on $B_0$ as $\lambda \to \infty$,
 or for a  sample $\U_{n,1}$  
of non-random size $n$ from the
uniform distribution on $B_0$ as $n \to \infty$.
 Here, we just state a  result of the second type.
\begin{theo}
\label{Xncoro}
Suppose $B_0 \in \RR(\R^d)$ with $|B_0|>0$. 
Set $\la_0 := |B_0|^{-1}$.
 Suppose $\HH^1, \ldots, \HH^k$ are 
 point process set functions,
satisfying the strong stabilization conditions
$(\ref{sstab1})$, $(\ref{sstab2})$, the  moments condition
 $(\ref{ubdmom})$, the uniform bound
$(\ref{polybd})$, and  homogeneity of order $\gamma$
$(\ref{eqhomog})$ for some $\gamma \in \R$.
 Suppose that $A_1,\ldots, A_k$ are 
sets in $\RR(B_0)$, and let $\Tau^{\la_0} =(\tau_{ij}^{\la_0})_{i,j=1}^k$
 be given by $(\ref{taudef})$.
Then for $i,j \in \{1, \ldots, k\}$
 \bea
\limn n^{(2 \gamma/d)-1} \cov(\HH^i(\U_{n,1},A_i), \HH^j(\U_{n,1}, A_j)) =
\tau_{ij}^{\la_0}
\label{0724a}
\eea
 and as $n \to \infty$,
\bea
n^{(\gamma/d) -1/2} (\HH^i(\U_{n,1},A_i) - \E \HH^i(\U_{n,1},A_i)
)_{i=1}^k \tod 
\NN(0, \Tau^{\la_0}).
\label{0724b}
\eea
\end{theo}
Theorem  \ref{Xncoro} is easily proved by
applying Theorem \ref{2Hclt} with $t_n= n^{1/d}$,
 and using homogeneity of $\HH^i $ to deduce that
$$
( n^{\gamma/d} h^i(\U_{n,1},A_i) )_{i=1}^k
\eqd 
(h^i ( \U_{n,n^{1/d}} , n^{1/d} A_i) )_{i=1}^k.  
$$

Many
 applications 
are concerned with functionals of graphs
 of the form $G:=G(\X)$ defined for each locally finite point set
 $\X \subset \R^d$ (a locally finite subset of $\R^d$ is  one with no limit
 point), where $G(\X)$ has vertex set $\X$.
See Sections \ref{secmst} and \ref{secnngs} for examples.

 We shall say $G$ is {\em translation invariant} if
translation by $y$ is a graph isomorphism from
$G(\X)$ to  $G(\tau_y(\X))$ for all $y\in \R^d$ and all
locally finite  point sets $\X$.
We shall say $G$ is {\em scale invariant} if
$G(a\X)$ is isomorphic to $G(\X)$ 
  for all  $\X$ and all $a >0$.

We use the following notion of stabilization for these graphs.
Given $G$, and given a vertex $x \in \X$, let 
$\Ed^+(x;\X)$ be
the set of edges of $G(\X)$ which are not edges of $G(\X \setminus \{x\})$, and
let
$\Ed^-(x;\X)$ be
the set of edges of
$G(\X \setminus \{x\})$
 which are not edges of 
 $G(\X)$.
Let $\Po_\la^0 := \Po_\la \cup \{\0\}$.
Our stabilization condition for graphs is
that there exists an almost
 surely finite random variable $R$ such that
\bea
\Ed^+(\0;\Po^0_\la  ) = \Ed^+ (\0; (\Po_\la^0 \cap B_R(\0) )
\cup \A)
\label{gstab1}
 \eea
and 
\bea
\Ed^-(\0;\Po_\la^0  ) = \Ed^-(\0; (\Po_\la^0 \cap B_R(\0) )
\cup \A)
\label{gstab2}
 \eea
for all finite $\A \subset \R^d \setminus B_R(\0)$.

The stabilization conditions (\ref{gstab1}), (\ref{gstab2})
say that
the local behavior
 of the graph in a bounded region   is unaffected by points
beyond a finite (though possibly random) distance from that region. As we
shall see,  the minimal spanning tree
(with the definition suitably extended from finite to locally finite
point sets) and the $k$-nearest neighbours,
and sphere of influence graphs
 all satisfy (\ref{gstab1}) and (\ref{gstab2}).

Another technical condition that turns out
to be  relevant to stabilization is {\em uniqueness of the infinite
component} for $G(\Po_\la)$ and for $G(\Po_\la^0)$.
  For a locally  finite point set
$\X$,
we say that uniqueness of the infinite component holds
for $G(\X)$ if there is almost surely at most a single
infinite component of $G(\X)$.

Given $G$,
we consider three types of functional based on $G$. 
Firstly, we consider the number of
components of $G(\X)$ 
with at least one vertex in $A$, which we denote $K^G(\X,A)$.

Second, functionals such as  total length of edges 
in $A$, number of edges in $A$, or
 number of edges in $A$ of  less than some specified length
may be interpreted as a
 total of $\phi$-weighted edge lengths in $A$, i.e.,
as a sum
 \bea
 L_\phi^G(\X,A):= \frac{1}{2} \sum_{x \in \X \cap A}
 \sum_{e =\{x,y\} \in G(\X)} \phi(|e|),
\label{phiw1} \eea
 for some appropriately specified function
 $\phi: (0,\infty)\to \R$.

Third,  we consider functionals such as the number of vertices
in $A$ of
some specified degree, or the number of components in $A$
with a specified number of vertices,
which are obtained by summing over
all vertices in $A$ some function of the local graph landscape
of $G$ (not the edge lengths) at that vertex.
To make this precise, let ${\cal K}$
denote the set of unlabelled connected rooted graphs
(i.e., connected graphs with a single vertex  distinguished and
denoted the {\rm root}).
For $\kappa \in \N$,
let ${\cal K}_{\kappa}$ denote the set of graphs in ${\cal K}$
which have all vertices a graph distance at most ${\kappa}$
from   the root (the graph distance between
two vertices is the minimal number of edges in
a path between them, or infinity if no such path exists).
 For $\kappa \in \N$,
let $B(\K_\ka)$ denote the class of all 
bounded real-valued functions  from ${\cal K}_\kappa$  to 
$\R$. For $\psi \in B(\K_\ka)$ and   for
 any vertex $x$ of any locally finite
point set  $\X \subset \R^d$,
let $G_{x,\ka}(\X)$ denote the rooted  subgraph of $G(\X)$ induced
by all vertices a graph distance at most $\kappa$ from $x$
(with root at $x$),
and let 
$$
V^G_\psi(\X,A) : = \sum_{x \in \X \cap A} 
\psi ( G_{x,\ka}(\X)).
$$

\begin{lemm}
\label{lem0516} Suppose $G$ is translation invariant and
satisfies the stabilization conditions $(\ref{gstab1})$ and $(\ref{gstab2})$. 
 Then if we set  $\HH(\X,A)= L^G_\phi(\X,A)$,
the stabilization conditions
$(\ref{sstab1})$ and $(\ref{sstab2})$ hold.
If instead, for some $\ka \in \N$ and $\ka \in B(\K_\ka)$, we set
 $h(\X,A) = V_\psi^G(\X,A)$ then, again,
conditions $(\ref{sstab1})$ and $(\ref{sstab2})$ hold.

Suppose in addition that
uniqueness of the infinite component holds for $G(\Po)$ and 
for $G(\Po^0)$;  then if we set
 $h(\X,A) = K^G(\X,A)$, then the
stabilization conditions $(\ref{sstab1})$ and $(\ref{sstab2})$ hold.
\end{lemm}
As we shall see in examples below, one can use Lemma
\ref{lem0516}
to check the applicability of Theorems \ref{1Hclt}, \ref{2Hclt} and 
\ref{Xncoro}
to a variety of point process functionals based on 
stabilizing graphs. 

\subsection{Marked point processes}
In the application in Section \ref{secOL},  we need
to consider the extension of the results
of the preceding section to functionals of {\em marked} point processes
with marks in the unit interval.
A {\em  marked point set} in $\R^d$ is a
 locally finite subset of $\R^d \times [0,1]$
 with no two elements having the same coordinate projection
onto $\R^d$.

If $\tilde{\X} = \{(x_i,t_i), i \geq 1\} \subset \R^d \times [0,1]$
 is a marked point set in $\R^d$, and $\X = \{x_i, i \geq 1\}$
is the corresponding unmarked set (i.e., the projection of
$\tilde{\X}$ onto $\R^d$), then we shall often abuse notation
slightly and write $\X$ for $\tilde{\X}$, keeping in mind
that each element $x_i$ of $\X$ carries a mark $t_i$. 
 Then all the notions
and results of the previous section carry through,
as we now describe.

For $y \in \R^d$,
the translation operator $\tau_y$ 
on marked point sets in $\R^d$ is to be understood
to preserve the values of all marks.  Then
the notion of a (marked) point process set function
$h(\X,A)$, defined for finite marked point sets
 $\X  $ in $\R^d$ and for $A \in \RR(\R^d)$,
is as given 
at the start of Section \ref{contsubsec}.
Also, the notion of translation invariance of
a graph $G(\X)$ is as defined in Section \ref{contsubsec}.
When we consider edge lengths and so on,
the vertex set of $G(\X)$ is still viewed as
a  subset of $\R^d$, not $\R^{d+1}$.
Also, it is to be understood that scalar multiplication
operator $\X \mapsto a\X$ on marked point sets in
$\R^d$, seen in the homogeneity condition (\ref{eqhomog}) 
for example,
leaves all marks unchanged.

In the marked setting,
the points of the $d$-dimensional point processes $\Po_{\la,t}$ and  
${\cal U}_{m,t}$ are to be understood
to carry marks which are each uniformly distributed
on $[0,1]$ and independent. 
Also, the inserted point at $\0$, when
defining add one
costs such as at (\ref{add1c}), is assumed to
carry an independent mark which is also uniformly 
distributed  on $[0,1]$.
The stabilization conditions (\ref{sstab1}) and (\ref{sstab2}) 
are to be understood to hold for any choice of 
values for the marks of points in ${\cal A}$. Likewise
the uniform bound (\ref{polybd}) is to be understood
to hold for any choice of the marks on $\X$.

With these interpretations,  
all of the results in Section \ref{contsubsec}
remain valid for marked point set functionals on
the marked point processes and stabilizing
graphs on the marked point processes.

\section{Applications of the general results}
\label{Secaps}
\allco
\subsection{Percolation}
Let $E=\{0,1\}$, let ${\cal E}$ be the power set of $E$
(i.e. the collection of all subsets of $E$),
with $P_0(\{1\})=p$ 
and $P_0(\{0\})=1-p$, $p\in (0,1)$ a fixed parameter. 
Let $X = (X_z)_{z \in \Z^d}$ and $X^\0$ (the same as $X$ but
with $X_0$ resampled) be as described in
Section \ref{secmainth} with this choice of $(E,{\cal E},P_0)$.
Let the sets $\Op$, $\Op'$ 
 (the random set of `occupied sites' induced by $X$ and by $X^\0$
respectively)  
be given by
$$
\Op := \{ z \in \Z^d : X_z=1\},~~~~~
\Op' : = \{ z \in \Z^d : X_z^\0=1\}
$$
(so that $\Op' \bigtriangleup \Op$ is either the empty set
or the set $\{\0\}$).

For any subset $S$ of $\Z^d$, let $G(S)$ be
the graph with vertex set 
$\{z \in S: X_z =1\}$, and
with edges between each pair of vertices at unit Euclidean
distance from each other. Then
 $G(\Op)$ is a Bernoulli site percolation process with parameter $p$
(bond percolation versions of
the results in this section also hold,
and are proved by similar means  taking $E= \{0,1\}^d$;
see \cite{Pclt}, page 1517). 

For background information on percolation see Grimmett \cite{Grimmett}.
Let $p_{\rm c}$ be the critical value of $p$, i.e.,  the
 supremum of the set of $p$ for which
the components of $G(\Op)$ are a.s. all finite.
Provided $d \geq 2$, it is known that $0 < p_\cc < 1$.

By the uniqueness of the infinite cluster in percolation
(see, e.g., \cite{Grimmett}), there is almost surely  
at most a single infinite component  
of $G(\Op)$. For later use, we denote the vertex set of
this infinite component of $G(\Op)$ by 
 $C_\infty$ (possibly the empty set), and
 denote the vertex set of
this infinite component of $G(\Op')$ by 
 $C'_\infty$.

Also for later use, observe that for $y,z \in \Z^d$,
\bea
(\tau_yX)_z =1 \Longleftrightarrow y+ z \in \Op \Longleftrightarrow
z \in \tau_{-y}(\Op).
\label{0903}
\eea

We shall give two applications of
Theorem \ref{basictheo} to percolation.
Suppose $B_0 \in \RR(\R^d) $ with $|B_0|>0$.
The next result adds to previously known central
limit theorems for 
the total number of components in $tB_0$ 
(see \cite{CG,Grimmett,Pclt,Z}),
and says that the number of components 
of $G(\Op \cap t B_0)$ 
in disjoint subregions of a large region $tB_0$ are asymptotically
normal and asymptotically independent of each other. It
is  of particular interest in the case when $p = p_\cc$ since
in this case, correlations  are not believed to decay exponentially.

For percolation and also for some of the other spatial graphs
that we consider, there are several ways to count
the `number of components' in a subregion $A$ of $\R^d$,
since one has to decide whether to include
components that lie only partially in $A$. In results
given here, such components are counted fully,
but the same results should hold if they were
counted only partially, or not at all.

\begin{theo} 
\label{perctheo}
Suppose $B_0 \in \RR(\R^d) $ with $|B_0|>0$.
For $t \geq 1$, $A \in \RR(B_0)$,
let $H_t(A)$ be the number of components of $G(\Op \cap tB_0)$
which include at least one vertex in $tA$.


Let $\Delta^H_\infty$ be the number of components of $G(\Op)$
that include at least one vertex at or adjacent to the origin,
minus the number of components of $G(\Op')$
that include at least one vertex at or adjacent to the origin.

 Then the conditions $(\ref{0414a})$, $(\ref{0414b})$ and $(\ref{4moments})$ for
Theorem $\ref{basictheo}$ are satisfied, and therefore the conclusions
$(\ref{limvar})$ and $(\ref{limdist})$ of that result are valid
$($with $H^i=H$ for all $i.)$
\end{theo} 

\noindent
{\bf Remark.} Following the approach of Cox and Grimmett \cite{CG} to
the central limit theorem for the number of components in $tB_0$,
 one could generalize Theorem \ref{perctheo} by taking $H_t(A)$ to
be of the form $\sum_{C \in {\cal C}(t,A)} \psi(C)$. Here
${\cal C}(t,A)$ denotes the 
set of $C \subset \Z^d$ such that $C$ is the vertex set of a 
component of $G(\Op \cap tB_0)$ which 
has at least one vertex in $tA$, and $\psi$ is some function defined
on finite $S \subset \Z^d$  such that $G(S)$ is  connected 
(in Theorem \ref{perctheo} we consider
the special case where $\psi$ is identically 1).
In this more general setting, by a modification
of the proof of Theorem \ref{perctheo}, one can still check the conditions
(\ref{0414a}), (\ref{0414b}), and (\ref{4moments}), and
hence apply Theorem \ref{basictheo} 
provided $\psi$ satisfies the following conditions:
  \begin{enumerate}
\item
 $\psi$ is translation-invariant, i.e. $\psi(\tau_y(S))=\psi(S)$
for all  $y \in \Z^d$ and all
$S \subset \Z^d$ such that $G(S)$ is connected.
\item
$\psi(S)$ converges to a finite limit as $|S| \to \infty$. 
\end{enumerate}
The above conditions imply that $\psi$ is bounded.
Unlike in \cite{CG}, we do not require
$\psi$ to be monotone here and we can take any $p \in (0,1)$,
including $p=p_c$. On the other hand, the corresponding
set of conditions on $\psi $ in  \cite{CG} does not include 
translation-invariance. \\

\noindent
{\em Proof of Theorem \ref{perctheo}.}
Let $t \geq 1$, $y \in \Z^d$, and $A \in \RR(B_0)$.
By (\ref{0903}), $H_{t,y}(A)$ is the number of components of
$G((\tau_{-y}\Op)  \cap tB_0)$ which intersect $tA$ (i.e., 
contain at least one vertex in $tA$). Hence, $H_{t,y}(A)$
is the number of components of
$G(\Op \cap \tau_y(t B_0))$ 
which intersect  $\tau_y(tA)$.

Thus  $- \Delta_{t,y}^H(A)$ is the increment in the
the number of components of $G(\Op \cap \tau_y(t B_0))$ 
which intersect  $\tau_y(tA)$
 when we resample $X_\0$ (i.e., when we replace the process $X$ by $X^\0$).

With  $\Delta^H_\infty$ defined in the statement of the theorem,
we assert that (\ref{0414a}) and (\ref{0414b}) hold. 
To verify (\ref{0414a}), suppose that
 $\liminf_{n \to \infty} (\tau_{y_n}(t_n A)) =\R^d$. 
Suppose first that $X_0=0$ and $X_*=1$. Then
there exists a (random) $N_1$ such that for   $n \geq N_1$, 
every pair of vertices lying adjacent to the origin and in the same
component of $G(\Op)$, is connected by an  path in
$G(\Op)$, all of whose  vertices lie  in $\tau_{y_n}(t_nA)$.
Then for all  $n \geq N_1$,
$\Delta^H_{t_n,y_n}(A) = \Delta^H_\infty$ as described above.
A similar argument applies in the case with
 $X_0=1$ and $X_*=0$, and for other cases
clearly  $H_{t_n,y_n}(A)=0$ for all $n$.
 Thus (\ref{0414a}) holds.

Next, suppose $\liminf(\tau_{t_n,y_n}(B_0\setminus A) )= \R^d$.
Suppose $X_0=0$.
There exists
a random $N_2$ such that for all 
large enough  $n \geq N_2$, the set 
 $\tau_{t_n,y_n}(B_0\setminus A) $
contains all finite components of $G(\Op)$
lying adjacent to the origin.

\begin{figure}[htbp]
\begin{center}
\includegraphics[width=8cm]{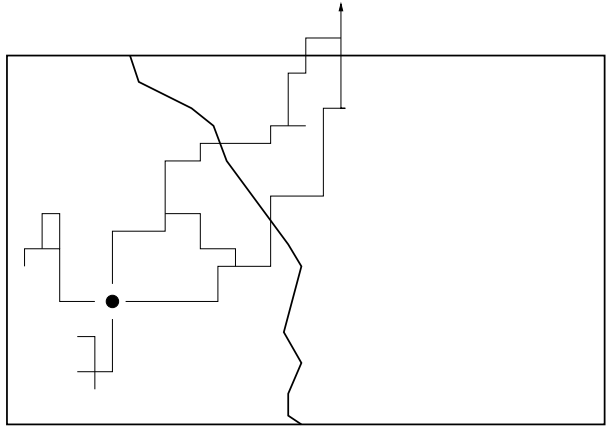}
\end{center}
\caption{The rectangle  represents the region $\tau_{y_n}(t_n B_0)$, 
and the region to the right of the jagged line is $\tau_{y_n}(t_nA)$.
The other rectilinear lines indicate the components
of $G(\Op)$ adjacent to the origin.}
\label{fig1}
\end{figure}

There exists a random $ N_3$ such that
 for all $n \geq N_3$, each pair of vertices of $C_\infty$ 
  which lie adjacent to the origin
is connected by a path in $G(\Op)$  all of whose vertices 
lie in the set 
 $\tau_{y_n}(t_n(B_0\setminus A) )$.
We assert that if
$n \geq \max(N_2,N_3)$,  changing of the 
value of $X_{\0}$ from $0$ to 1 
does not affect the number of components
of $G(\Op \cap \tau_{y_n}(t_nB_0))$
 that intersect $\tau_{y_n}(t_nA)$. This is because for
$n$ this big, any two occupied vertices adjacent to the origin
which are both connected by   paths in $G(\Op)$ to vertices in
 $\tau_{y_n}(t_n A)$, must be part of $C_\infty$
and therefore are connected by a path which avoids the origin,
so that they are already part of the same component
of $G(\Op \cap \tau_{y_n}(t_nB_0))$ even before we add a vertex at the
origin to $\Op$. In other words, 
$G(\Op \cap \tau_{y_n}(t_n B_0))$ has at most a single
component which intersects both the set $\tau_{y_n}(t_nA)$ and
and  the set of sites adjacent to the origin;
  see Figure \ref{fig1}.  The assertion  follows, and one argues similarly
for $X_*=0$. Thus (\ref{0414b}) holds. 

Since all vertices in $G(\Z^d)$ have degree $2d$,
 the absolute value of $\Delta_{t,y}(A)$ is uniformly bounded by
$2d-1$, and therefore the moments  condition (\ref{4moments})
is valid for any finite $\gamma$. Therefore  Theorem \ref{basictheo} is
applicable here. $\qed$ \\

We now consider the {\bf largest component}
of $G(\Op \cap tB_0)$, adding to the central limit theorem
for the largest component size given in \cite{Pclt}. A largest component 
of $G(\Op \cap tB_0)$ is a component such that no other 
component has more vertices.
There could be more than one largest
component; in the sequel, the ``vertices which lie
in a largest component'' means the vertices lying in the  union
of all largest components as defined above, while any discussion
of properties of ``{\bf the} largest component'' refers to the case
where there is a unique largest component.
 The following result says that the
distribution of the vertices lying in a largest component in $tB_0$ 
is asymptotically a white noise distribution. 
In it, we assume $B_0$ is rectangular,
i.e., that $B_0$ is a product of bounded intervals. Presumably,
the proof can be extended to other shapes of $B_0$.

\begin{theo}
Suppose that $p> p_\cc$, and that $B_0$ is rectangular.
For $t \geq 1$ and $A \in \RR(B_0)$,
let $H_t(A)$ be the number of vertices
of $tA$ which lie in a largest component of $G(\Op \cap tB_0)$. 
 Set
$$
\Delta^H_\infty = \left\{ \begin{array}{llr} 
\card(  C_\infty \setminus C'_\infty ) & {\rm if }& X_\0=1, X_*=0,
\\
-\card(  C'_\infty \setminus C_\infty ) & {\rm if }& X_\0=0, X_*=1,
 \\
0 & {\rm if } & X_\0 = X_*
\end{array}
\right.
$$
Then
 the conditions
 $(\ref{0414a})$, $(\ref{0414b})$ and $(\ref{4moments})$ for
Theorem $\ref{basictheo}$ are satisfied, and therefore the conclusions
$(\ref{limvar})$ and $(\ref{limdist})$ of that result are valid
$($with all $H^i=H.)$
\end{theo}

{\em Proof.}
Observe first that $\Delta^H_\infty$ is indeed almost surely finite.
For example, if
$X_\0 =1$ and $X_* =0$, then
 $C'_\infty \subseteq C_\infty$ and 
$C_\infty \setminus C'_\infty$
consists
(in the case where $\0 \in C_\infty$)
 of those finite components of $G(\Op')$ lying adjacent to
the origin, along with the origin itself (with $C_\infty \setminus
C'_\infty = \emptyset $ in the case where $\0 \notin C_\infty$).

Observe also that for any $t \geq 1$, $y \in \Z^d$ and
$A \in \RR(\B_0)$, by (\ref{0903}), $H_{t,y}(A)$ is the number of
vertices in $tA$ in a largest component of $G(\tau_{-y}(\Op)\cap tB_0)$,
and so is  the number of vertices in 
 $\tau_y(tA)$ in a largest  component of
$G(\Op \cap \tau_y(tB_0))$.

In what follows, a few plausible (and actually true)
 facts are stated without proof. For details
of their proofs, see \cite{Pclt}.

Suppose  that $\liminf_{n \to \infty} (\tau_{y_n}(t_n B_0)) = \R^d$.
Then with probability tending to 1, the largest component
of $G(\Op \cap \tau_{y_n}(t_nB_0))$ is unique, 
and is the largest  component of
$G(C_\infty \cap \tau_{y_n}(t_nB_0))$, and if $\0 \in C_\infty$ 
then the largest component of $G(\Op \cap \tau_{y_n}(t_nB_0))$ is the
component of 
$G(C_\infty \cap \tau_{y_n}(t_nB_0))$ containing the origin.

Suppose that  $\liminf_{n \to \infty} (\tau_{y_n}(t_n A)) = \R^d$.
Then with probability  1, the set $C_\infty \bigtriangleup C'_\infty$
is  contained in $\tau_{y_n}(t_nA)$ for all large enough $n$.
Hence, the probability that  $\Delta_{t_n,y_n}(A)$ is equal to
$\Delta^H_\infty$ defined above tends to 1, and so
(\ref{0414a}) holds.

Suppose  that
 $\liminf_{n \to \infty} (\tau_{y_n}(t_n (B_0 \setminus A))) = \R^d$.
Then with probability  1, the set $C_\infty \bigtriangleup C'_\infty$
is  contained in $\tau_{y_n}(t_n(B_0 \setminus A))$ for
 all large enough $n$.
If 
 $C_\infty \bigtriangleup C'_\infty \subseteq
\tau_{y_n}(t_n(B_0 \setminus A))$ and also
the largest component of $G(\Op \cap \tau_{y_n} (t_n B_0))$ is the
sole component of $G(C_\infty \cap \tau_{y_n} (t_nB_0))$ containing 
a vertex adjacent to the origin,
and also
the largest component of $G(\Op' \cap \tau_{y_n} (t_n B_0))$ is the
sole component of $G(C'_\infty \cap \tau_{y_n} (t_nB_0))$ containing 
a vertex adjacent to the origin,
 then 
changing the value of $X_\0$ from 1 to 0 will not remove any vertices
of the largest component lying in $\tau_{y_n}(t_nA)$ so that 
it does not change the value of $H_{t_n,y_n}(A)$. 
Hence the probability that  $\Delta^H_{t_n,y_n}(A)$ is equal to
zero tends to 1, and so (\ref{0414b}) holds.

We need to check the moments condition (\ref{4moments}).
Most of the ingredients in the proof of this are given in
 the proof of Theorem 3.2 of 
\cite{Pclt}. The main  difference is that we now need to account 
for a possible {\em decrease} in the number of
elements in $\tau_y(tA)$ of a largest component when
we change the status of site $\0$ from `vacant' to `occupied'
(in \cite{Pclt} we needed only to
 consider the  largest component size, which by contrast
really is monotone in $X_\0$).
Such a decrease could happen
either if 
 $G(\Op \cap \tau_y(tB_0) \setminus \{\0\} )$ 
has more than one largest component, or
 if two or more components 
of $G(\Op \cap \tau_y(tB_0) \setminus \{\0\} )$ 
lying adjacent to $\0$,
when merged, form a component larger than and disjoint from the previous
largest component.
 However, the probability of either of these possibilities occurring
decays exponentially in $t^{d-1}$ (see, e.g., Theorem 4 in Penrose and
Pisztora \cite{PP}), and using this we can check (\ref{4moments})
here. $\qed$

\subsection{The minimal spanning tree}
\label{secmst}
Given a locally finite set $\X \subset \R^d, \ d \geq 2$, and given $a >0$, let
${\cal G}_a(\X)$  be the graph with vertex set $\X$ and with edge set
$\{\{x,y\}: |x-y| < a\}$.  Let $\MST(\X)$ be the graph with
vertex set $\X$ obtained by including each edge
$\{x,y\}$ such that $x$ and $y$ lie in different components of
 ${\cal G}_{|x-y|}(\X)$ and
at least one of these components is finite.
 If $\X$ is finite with
distinct inter-point distances, then $\MST(\X)$  is the
  {\em minimal spanning tree} on $\X$, i.e. the
connected graph with vertex set $\X$ of minimal total edge length;
see Aldous and Steele (\cite{AS}, Lemma 12).
Clearly $\MST(\X)$ is translation and scale invariant.

Recall the definitions of $V_\psi^G$, $L_\phi^G$, and
 $B(\K_\ka)$  from Section \ref{contsubsec}.
The first part of the following result
tells us that  the totals of a local graph landscape 
function (for example, the numbers of leaves), 
summed over points of the random minimal spanning tree in
disjoint regions, scaled and centred, are asymptotically
independent normals.
The second part says that the totals of $\phi$-weighted edges of the 
random minimal spanning tree in disjoint regions,
scaled and centred, are asymptotically independent normals.
In this result, 
say $\phi$ is {\em polynomially bounded}
if there exists a constant $c$ such that
$|\phi(r)| \leq c(1+r)^c$ for all $r>0$.

\begin{theo}
\label{theomst}
Suppose $G(\X)$ is $\MST(\X)$.
Let $\la >0$, and
suppose $B_0 \in \RR(\R^d)$ with $|B_0|>0$.

Let $\kappa \in \N$, and suppose $\psi \in B(\K_\ka)$.
If we set
 $h(\X)= V_\psi^G(\X,A)$, then $h$ satisfies all the
conditions 
$(\ref{sstab1})$, $(\ref{sstab2})$, $(\ref{ubdmom})$,
$(\ref{polybd})$, and $(\ref{eqhomog})$ $($with $\gamma=0)$
of Theorems $\ref{1Hclt}$, $\ref{2Hclt}$,
and $\ref{Xncoro}$ and therefore satisfies their 
  conclusions $(\ref{0516a})$, $(\ref{0516b})$, $(\ref{varlim})$
$(\ref{2clteq})$,
$(\ref{0724a})$ and $(\ref{0724b})$  
 $($with $h^j(\X,A) =V_\psi^G(\X,A)$ for all $j$ and with $\gamma=0)$.

 Suppose instead that we set $h(\X,A) = L_\phi^G(\X,A)$ for
some $\phi: (0,\infty) \to \R$. Then
the stabilization conditions $(\ref{sstab1})$  and $(\ref{sstab2})$ hold.
If  $\phi$ is bounded, or if 
$\phi$ is polynomially bounded and $B_0$ is convex,
then the moments condition $(\ref{ubdmom})$ holds and so
  Theorems $\ref{1Hclt}$ and $\ref{2Hclt}$ apply and their
  conclusions $(\ref{0516a})$, $(\ref{0516b})$, $(\ref{varlim})$ and
$(\ref{2clteq})$
 $($with $h^j =h$ for all $j)$ hold.

If also $\phi(r)= r^\alpha$ for some constant $\alpha$, 
then the homogeneity hypothesis $(\ref{eqhomog})$ holds and hence the
 conclusions $(\ref{0724a})$ and $(\ref{0724b})$ of Theorem
$\ref{Xncoro}$ are valid with $\gamma =\alpha$.
\end{theo}
The proof of this uses the following lemma which we shall
use again later.
\begin{lemm}
\label{lemregular}
If $B_0 \subseteq \R^d $ is bounded and convex with $|B_0|>0$, then
\bea
\inf_{x \in B_0, r \in (0,1]} 
r^{-d} |B_r(x) \cap B_0| >0.
\label{0903a}
\eea
\end{lemm}
{\em Proof.}
The assumptions on  $B_0$ imply that
$B_0$ has non-empty interior, so that there exists
 a ball $B$ contained in $B_0$.
For any $x\in B_0$, the convex hull of $\{x\} \cup B$ is
contained in $B_0$, and since $B_0$ is bounded
the angle subtended by this cone-like set at $x$
is bounded away from zero, so the result follows.
$\qed$\\

\noindent
{\em Proof of Theorem \ref{theomst}.}
Condition (\ref{gstab1}) follows from Lemma 2.1 of \cite{PY4}.
Condition (\ref{gstab2}) is more complicated but follows from
the proof of Proposition 1 of Lee \cite{LeeII}.
Therefore, Lemma \ref{lem0516} of the present paper
can be applied to give us the conditions
(\ref{sstab1}), (\ref{sstab2}) in the case where either
$h(\X,A)= V_\psi^G(\X,A)$ or  
$h(\X,A)= L_\phi^G(\X,A)$. 

Given a finite set $\X \subset \R^d$,
consider the effect on the minimal spanning tree MST$(\X)$ 
of adding a point at the origin $\0$.  Let edges of MST$(\X \cup \{\0\})$
that are not in MST$(\X)$ be denoted {\em added edges}, and let edges of 
MST$(\X)$ that are not in MST$(\X \cup \{\0\})$
be denoted {\em deleted edges}.

By the revised add and delete algorithm of Lee \cite{LeeI}, 
the  added edges are precisely those incident to
$\0$ in MST$(\X \cup \{\0\})$, and 
there are fewer deleted edges than added edges.
Moreover, there is a uniform non-random  bound on vertex
degrees in the minimal spanning tree
 (see \cite{AS}), and hence there is a uniform bound both on
the number of added edges and on the number of deleted edges.
The moments condition 
(\ref{ubdmom}) for $h(\X,A)= V_\psi^G(\X,A)$ 
is immediate
 from these remarks. 
Moreover,
if $\phi$ is bounded, then 
(\ref{ubdmom}) for $h(\X,A)= L_\phi^G(\X,A)$ 
also follows from these remarks. 
 
Suppose that $\phi$ is polynomially bounded and
$B_0$ is convex.
By the preceding remarks,
to prove (\ref{ubdmom}) in this case, it suffices
 to show that for any $K>0$, there is a 
deterministic uniform bound on the
 $K$th moment of the length of the longest added edge
when a point at $\0$ is inserted into $\U_{m,t}$ with $t \geq 1$,
$\0 \in \tau_x(tB_0)$ and $m/(\lambda t^d |B_0|)$ in the range $[1/2,3/2]$,
and likewise
for the the longest deleted edge.

We assert that the longest deleted edge in a finite set $\X$ is
at most twice as long as the longest added edge. To see this,
suppose that
 $\{X,Y\}$ is a deleted edge.
Then, since all added edges are incident to the added point at $\0$, 
and there must be a path from $X$ to $Y$ in MST$(\X\cup \{\0\})$,
there exist points $X',Y'$ in $\X$  such that $X',Y'$
are both adjacent to $\0$ in $\MST(\X \cup \{\0\})$, and such
 that there  a path in $\MST(\X)$ from $X$ to $X'$, and a path 
in $\MST(\X)$ from $Y$ to $Y'$. 
By the triangle inequality,
$|X'-Y'|$ is at most twice the
length of the longest added edge, and also
$|X-Y| \leq |X'-Y'|$ since otherwise we could start
with $\MST(\X)$, then replace edge $\{X,Y\}$ by $\{X',Y'\}$
to obtain a spanning tree on $\X$ of smaller total length, a
contradiction. This completes the proof of the assertion.

Thus, to check (\ref{ubdmom}) for $h(\X,A)= L_\phi^G(\X,A)$ 
when $\phi$ is polynomially bounded and $B_0$ is convex, 
 it suffices to prove the $K$th moment of the longest added edge
is uniformly bounded. This can be proved by an argument along the lines
of Lemma 2.1 of Yukich \cite{Y} (using convexity of $B_0$, and
(\ref{0903a})).

The uniform bound (\ref{polybd}) is trivial for $V_\psi^G$ and 
also holds for $L_\phi^G$ since $\phi$ is assumed polynomially bounded. 
The homogeneity condition (\ref{eqhomog}),  for 
$V^G_\psi$ or $L^G_\phi$  with $\phi(r)= r^\alpha$, follows
from the fact that the graph $\MST(\X)$ is scale invariant.  $\qed$\\

{\em Remark.} In the univariate case  Kesten and
 Lee (\cite{KL}, Theorem 2) give a CLT for $L_\phi^G(B_0)$
when $\phi$ is monotonically increasing. Our result
gives central limit theorems, for this univariate
case, without this restriction on $\phi$.\\

{\em Remark.}  Consider the empirical distribution of
MST edge lengths. That is, let $N_s(\X)$ be the number
of edges  of $\MST(\X)$ of length less than $s$. 
For arbitrary $B_0$, 
the finite-dimensional distributions
 of the scaled, centred empirical process
\bea
(t^{-d/2}(N_s(\Po_{\la,t}) - \E N_s(\Po_{\la,t})), ~~~~ s > 0)
\label{1121a}
\eea 
converge as $t \to \infty$ to those of a Gaussian process.
Moreover,
the finite-dimensional distributions
 of the corresponding 
scaled, centred empirical process
for a binomial sample, namely
\bea
(n^{-1/2}(N_{n^{-1/d}s}(\U_{n,1}) - \E N_{n^{-1/d} s}(\U_{n,1})), ~~~~ s > 0)
\label{1121b}
\eea 
converge as $n \to \infty$ to those of another Gaussian process.

 In the case of the first empirical process (\ref{1121a}),
this follows by taking an arbitrary set of positive `times'
$s_1, \ldots, s_k$ and applying Theorem \ref{1Hclt} with 
$h^i(\X,A)= L_{\phi^i}^G(\X,A)$ and $\phi^i(r) = {\bf 1}_{\{r \leq s_i\}}$
for $(1 \leq i \leq k)$.
 The limiting Gaussian process in this case
has covariance function given by the function
$\E [Y_sY_{s'}], s,s' > 0$,
where we set
\bea
Y_s = \E [ \delta_s(\infty)|\F],
\label{1121c}
\eea
 with $\delta_s(\infty)$ denoting the
stabilizing limit for the functional  
$h(\X,A)= L_{\phi_s}^G(\X,A)$ with $\phi_s(r): = {\bf 1}_{\{r \leq s\}}$,
and with $\F$ denoting the $\sigma$-field generated by
the Poisson configuration in the half-space 
$\{ x = (x^1,\ldots,x^d) \in \R^d : x_1 < 0\}$.

In the case of the  second empirical process
(\ref{1121b}), the convergence of finite-dimensional
distributions follows similarly but this time using
Theorem \ref{2Hclt} and the observation that
$(N_s(\U_{n,n^{1/d}}), s >0) $ has the
same distribution as
$( N_{n^{-1/d} s}(\U_{n,1}), s >0)$.
The limiting covariance  function
is this time given by $\Cov(Y_s,Y_{s'}), s, s' >0$, with $Y_s$
given once more by (\ref{1121c}), but now with 
$\lambda = |B_0|^{-1}$.

\subsection{Nearest-neighbour type graphs}
\label{secnngs}

Let $k \in \N$.
The $k$-nearest neighbour graph
($k$-NNG) on  a locally finite set 
$\X \subset \R^d$  is obtained  by  including an undirected edge
connecting each vertex $x \in \X$ to each of its $k$ nearest
neighbours (using the lexicographic ordering as a tie-breaker
in the event of a tie). We also consider the sphere of
influence graph (SIG), in which, denoting 
the distance from $x \in \X$ to its nearest neighbour by
$R_x$,
 we connect vertices $x,y$ of $\X$ 
 by an edge if and only if
 $B_{R_x}(x) \cap B_{R_y}(y) \neq \emptyset $.

White noise limits for 
functionals such  as $L_\phi^G$ and $V_\psi^G$
(defined in Section \ref{contsubsec}) can
be derived using either the results in this paper or by
other methods based on exponential decay,  as in
\cite{BY}. We concentrate here on the component
count $K^G(\X,A)$, for which exponential decay is not so clear. 

\begin{theo}
Suppose $G(\X)$ is defined  to be the $k$-NNG on $\X$. 
Let $\la >0$ and let $B_0 \in \RR(\R^d)$ with $|B_0|>0$.
Then
$h(\X,A) = K^G (\X,A)$
 satisfies all the conditions
$(\ref{sstab1})$, $(\ref{sstab2})$, $(\ref{ubdmom})$,
$(\ref{polybd})$, and $(\ref{eqhomog})$ $($with $\gamma=0)$
for Theorems $\ref{1Hclt}$, $\ref{2Hclt}$, and
$\ref{Xncoro}$.
  Therefore, $h(\X,A) = K^G (\X,A)$ satisfies the conclusions
$(\ref{0516a})$, $(\ref{0516b})$, $(\ref{varlim})$, $(\ref{2clteq})$,
$(\ref{0724a})$, and $(\ref{0724b})$
 in those  results $($with $h^i=h$ for all $i$, and with $\gamma=0)$.
\label{th0516a}
\end{theo}
{\em Proof.}  
There is a deterministic uniform upper bound on the degree of
vertices in $G(\X)$ (see, e.g., Lemma 8.4 of  Yukich \cite{Ybk}),
and hence  a uniform deterministic bound
on the change in the number of components of $G(\X)$ caused by inserting
a single point; the 
 moments condition (\ref{ubdmom}) follows.

We can obtain the stabilization conditions (\ref{sstab1}) and (\ref{sstab2})
by using Lemma \ref{lem0516}. This result is  applicable because 
conditions (\ref{gstab1}) and (\ref{gstab2}) hold by
the proof of Lemma 6.1 of \cite{PY}, while
uniqueness of the infinite
component in $G(\Po)$ and $G(\Po^0)$ holds by Lemma 6.4 of 
\cite{PY}. The uniform bound (\ref{polybd})
and the homogeneity (\ref{eqhomog}) are obvious.
$\qed$ \\

\begin{theo}
\label{th0516b}
Let $\la >0$, and suppose $B_0\subset \R^d$ is convex and bounded
with $|B_0|>0$.
Suppose $G(\X)$ is defined  to be the SIG on $\X$.  Then
 $h(\X,A) = K^G (\X,A)$
 satisfies all the conditions 
$(\ref{sstab1})$, $(\ref{sstab2})$, $(\ref{ubdmom})$,
$(\ref{polybd})$, and $(\ref{eqhomog})$ $($with $\gamma=0)$
for Theorems $\ref{1Hclt}$,   $\ref{2Hclt}$, and
 $\ref{Xncoro}$. 
Therefore, $h(\X,A) = K^G (\X,A)$ satisfies the conclusions
$(\ref{0516a})$, $(\ref{0516b})$, $(\ref{varlim})$, $(\ref{2clteq})$,
$(\ref{0724a})$, and $(\ref{0724b})$
 in those  results $($with $h^i=h$ for all $i$, and with $\gamma=0$$)$.
\end{theo}
{\em Proof.}  
We need to check the moments condition (\ref{ubdmom}). 
This can be done as in the proof of Theorem 7.4 of \cite{PY}.
The `regularity' condition in that result is implied 
by the condition (\ref{0903a})  here;
see the remarks on page 1010 of \cite{PY}.  

The rest of the argument is similar to that for 
Theorem \ref{th0516a}. 
For conditions (\ref{gstab1})
and (\ref{gstab2}), see the proof of Lemma 7.1 of \cite{PY}.
 For uniqueness of the infinite
component of $G(\Po)$, see Theorem 7.3 of 
\cite{PY}.  The infinite component of $G(\Po^0)$ 
is also unique, by the proof of Lemma 6.4 (b) of \cite{PY}.  $\qed$

\subsection{The on-line nearest neighbour graph}
\label{secOL}

Suppose $X_1,X_2,\ldots, X_n $ are
points in $\R^d$.
 In the 
{\em on-line nearest neighbour graph} (or on-line NNG for short), 
the points are assumed to arrive sequentially
 and each point
$X_i$, $i \geq 2$, is connected by an undirected edge
 to its nearest neighbour in
the set of preceding points in the sequence $\{X_1,\ldots, X_{i-1}\}$
(using the lexicographic ordering on $\R^d$ to break any ties).
The resulting graph is a tree, which we will denote
the {\em on-line nearest neighbour graph} on the
sequence of points $(X_1,X_2,\ldots,X_n)$.
 One could also consider the on-line $k$-nearest neighbour
graph defined analogously, with each point $X_i$ connected to its
$k$ nearest neighbours in $\{X_1,\ldots,X_{i-1}\}$ if $i >k$,
and  connected to each of $X_1,\ldots,X_{k-1}$ if $i \leq k$.
 In our case, the points in the sequence
 will be random, independent and uniformly distributed
over $B_0$ or over $tB_0$. 

The on-line nearest neighbour graph on random points 
is  a natural growth model for spatial random graphs,
although it was apparently introduced only recently, by
Berger et al. (\cite{BBBCR}, Section 3).
There, the motivation comes from the search for a simple
model of scale-free networks, the graph being itself
a simplification of a model of Fabrikant et al. 
\cite{FKP}.

To fit this graph into our present setup, consider a 
{\em marked} random finite point set $\X$ in $\R^d$, where 
each point $X$ of $\X$ carries a random mark $\TT_X$
 which is uniformly distributed on $[0,1]$,  independent of the other marks
and of the point process $\X$. The points are listed in
increasing order of mark, i.e., the marks represent time
of arrival. With this ordering,
we connect each point of $\X$ to the nearest point that
precedes it in the ordering, to obtain a graph which we also
call the on-line nearest neighbour graph on the marked point set $\X$.
This definition extends to  infinite but locally finite  
point sets.

Clearly the  on-line NNG
on ${\cal U}_{m,t}$ (defined via marked point processes)
has the same
 distribution as the on-line NNG (with the first definition)
on a sequence $X_1,\ldots, X_m$ of independent uniform
points on $tB_0$. Likewise, the on-line 
NNG for $\Po_{\la,t}$ (defined via marked point processes) has
the same distribution as the on-line NNG on
$(X_1,\ldots,X_N)$  with $\{X_i\}$ independent uniform
points on $tB_0$ and  $N$ an independent Poisson
variable with parameter $\lambda t^d |B_0|$. 

As it turns out, the on-line nearest neighbour graph
is  a nice example of our methods because it
is stabilizing but has only polynomially decaying
correlations (i.e., a polynomially decaying tail on
the distribution of the radius of stabilization).
 In the following discussion, although
we think of the graph as undirected, we shall refer
to an edge connecting marked points $X$ and $Y$
with $\TT_X > \TT_Y$, as an {\em outgoing edge}
from $X$ and an {\em incoming edge} to $Y$.
Each vertex (except one if $\X$ is finite) has
a single outgoing edge.

The existence of an almost surely finite
radius of stabilization satisfying (\ref{gstab1}) and
(\ref{gstab2}) will be shown later on.  To see that
its distribution does not have an exponentially
decaying tail,
 let $L$  be the length of the
outgoing edge from the  origin in the on-line nearest neighbour graph
on $\Po^\0$. With $\pi_d$ denoting the volume of the
unit ball in $\R^d$, we have
\bean
P[ L \geq r] = \int_0^1 \exp(-\lambda  \pi_d r^d t ) dt
= \lambda^{-1} \pi_d^{-1} r^{-d} \int_0^{\lambda \pi_d r^d}  e^{- u} du,
\eean
which shows that the tail of the distribution of
 $L $ decays only polynomially,
and $L$ is clearly a lower bound for any radius of
stabilization. 

Recall the definitions of $V_\psi^G$, $L_\phi^G$, and
 $B(\K_\ka)$  from Section \ref{contsubsec}.
The following result says that for certain $\phi$ the totals
 of $\phi$-weighted edges of the 
on-line nearest neighbour graph on random points in disjoint regions,
scaled and centred, are asymptotically independent normals, and
likewise for the totals of any bounded 
function of vertex degrees  summed over vertices in disjoint regions.

\begin{theo}
\label{thOL} Suppose $G(\X)$ is the on-line NNG  on $\X$.
Let $\la  > 0$ and
suppose $B_0 \subset \R^d$ is convex and bounded with $|B_0|>0$.
For  any
  $\psi \in B(\K_1)$,
 if we set
 $h(\X)= V_\psi^G(\X,A)$ then $h$ satisfies all the
conditions 
$(\ref{sstab1})$, $(\ref{sstab2})$, $(\ref{ubdmom})$,
$(\ref{polybd})$, and $(\ref{eqhomog})$ $($with $\gamma=0)$
of Theorems $\ref{1Hclt}$, $\ref{2Hclt}$,
and $\ref{Xncoro}$ and therefore satisfies their 
  conclusions $(\ref{0516a})$, $(\ref{0516b})$, $(\ref{varlim})$
$(\ref{2clteq})$,
$(\ref{0724a})$ and $(\ref{0724b})$  
 $($with $h^j(\X,A) =V_\psi^G(\X,A)$ for all $j$ and with $\gamma=0)$.

 Suppose instead that we set $h(\X,A) = L_\phi^G(\X,A)$ for
some $\phi:(0,\infty) \to \R$. Then
the stabilization conditions $(\ref{sstab1})$  and $(\ref{sstab2})$ hold.
If  $\phi$ satisfies the growth bound 
\bea
\sup_{ r >0}
\left( (1+r)^{-\alpha}  |\phi(r)|   \right)  < \infty, ~~~~{\rm some}
~ \alpha < d/4 ,
\label{growthbd}
\eea
then the moments condition
$(\ref{ubdmom})$ holds and so
  Theorems $\ref{1Hclt}$ and $\ref{2Hclt}$ apply and their
  conclusions $(\ref{0516a})$, $(\ref{0516b})$, $(\ref{varlim})$ and
$(\ref{2clteq})$
 $($with $h^j =h$ for all $j)$ hold.

If also $\phi(r)= r^\alpha$ for some constant $\alpha < d/4$, 
then the homogeneity hypothesis $(\ref{eqhomog})$ holds and hence the
 conclusions $(\ref{0724a})$ and $(\ref{0724b})$ of Theorem
$\ref{Xncoro}$ are valid with $\gamma = \alpha$. 
\end{theo}
{\bf Remarks.} 
Provided $B_0$ is convex,
 Theorem  \ref{thOL} gives, among other things,
a central limit theorem for the number of vertices of any
fixed degree in the on-line NNG on $\U_{n,1}$ or on
$\Po_{\la,t}$. 
Since  any bounded function of the edge lengths
satisfies the growth bound (\ref{growthbd}),
Theorem \ref{thOL} also enables us to
obtain similar functional central limit theorems 
results on the empirical
distributions  of edge lengths in the on-line 
nearest neighbour graph, to those described in the
preceding section for the minimal spanning tree.

Provided $d > 4$, Theorem \ref{thOL} gives us
a central limit theorem for the total
length of the on-line NNG on $\U_{n,1}$ or on
$\Po_{\la,t}$ (since the function $\phi(r)=r$ satisfies
the growth bound (\ref{growthbd})).
 This leaves open the question of the asymptotic
behaviour of the total length of the on-line NNG
on $\U_{n,1}$, in dimensions $d \leq 4$.  As mentioned earlier,
it is likely that 4th moments condition (\ref{ubdmom})
can be replaced by a $2+\eps$ moments condition
in Theorems \ref{1Hclt} and \ref{2Hclt}. If this can
be done, the total length of the
 on-line NNG on $\U_{n,1}$ will satisfy a central
limit theorem for $d=4$ or $d=3$. 
We suspect that a central limit theorem also holds
for $d=2$, but do not have  a proof. 
We believe that the limiting distribution of
the (centred) total length
is non-normal for $d=1$;  
Penrose and Wade \cite{PW} 
have shown this to be the case 
for a related graph in which 
$X_i$
is joined to its nearest neighbour to the {\em left}
in the set $\{X_1,\ldots,X_{i-1}\}$.

Theorem \ref{thOL}
also carries through to the on-line $k$-nearest neighbour
graph, although we give a proof only in the
case $k=1$. \\

The proof of Theorem \ref{thOL} uses the following three lemmas. 
The first two of these are purely geometric in nature.
Given distinct points $x,y \in \R^d$, let $C_{x,y}$ denote the cone 
with its point at $x$ and
with angular radius $\pi/12$, centred on the half-line from $x$
passing through $y$. 
\begin{lemm}
\label{lemgeo}
Suppose $B_0 \subset \R^d$ is convex and bounded with $|B_0|>0$. Then
$$
\inf_{x,y \in B_0: x \neq y} \frac{|C_{x,y} \cap B_{|y-x|}(x)
\cap B_0|}{|y-x|^d} >0.
$$ 
\end{lemm}
{\em Proof.} Take $x,y \in B_0$ with $|y-x|=r>0$.
Let $z = (x+y)/2$. By convexity $z \in B_0$, and 
geometrical considerations show that 
$$
B_{(r/2)\sin(\pi/12)}(z) \subseteq C_{x,y} \cap B_r(x),
$$
so the result follows from Lemma \ref{lemregular}.
$\qed$

\begin{lemm}
\label{lem0909}
Suppose $C$ is an open cone in $\R^d$, of angular
radius $\pi/6$, with its point at $x \in \R^d$. Then 
for $y \in C$ and $z \in C$ we have $|z-y|< \max(|z-x|,|y-x|)$.
\end{lemm}
{\em Proof.} Assume without loss of generality that
$x =\0$ and that $|y|\leq |z|$.
Let $\theta$ be the angle $z\0 y$. Then $\theta < \pi/3$, so 
$\cos \theta > 1/2$ and
by the cosine rule,
$$
|z-y|^2 = |z|^2 + |y|^2 - 2|y|\cdot|z| \cos \theta < |z|^2. ~~\qed \\
$$

Recall from Section \ref{contsubsec} that if a graph $G(\X)$
is defined for locally finite point sets $\X \subset \R^d$, then 
for $x \in \X$,
 ${\cal E}^+(x;\X)$ denotes the set of edges
of $G(\X )$ which are not
edges of $G(\X \setminus \{x\})$,
and ${\cal E}^-(x;\X)$ denotes the set of edges
of $G(\X \setminus \{x\})$ which are not edges of $G(\X )$.
Also, $|e|$ denotes the Euclidean length of edge $e$.
\begin{lemm}
\label{lemOL}
Let  $G(\X)$ denote the on-line nearest neighbour 
graph on the marked point set $\X$, and suppose $B_0 \subset \R^d$ is convex
and bounded with $|B_0|>0$.
 Let $\lambda >0$
and let $\mu_{\la,t} : = \lambda t^d |B_0|$. Let $0 \leq \alpha < d/4$. Then 
\bea
\sup_{t \geq 1, x \in tB_0}  
\sup_{m  \in [ \mu_{\la,t}/2, 3 \mu_{\la,t}/2] \cap \Z }  
\E \left[ \left(\sum_{ e \in {\cal E}^+(x; \U_{m,t} \cup \{x\} ) }
(1 + |e|)^{\alpha} \right)^4
\right]
 < \infty
\label{1128b}
\eea
and
\bea
\sup_{t \geq 1, x \in tB_0}  
\sup_{m  \in [ \mu_{\la,t}/2, 3 \mu_{\la,t}/2] \cap \Z}  
\E \left[ \left( \sum_{ e \in {\cal E}^-(x; \U_{m,t} \cup \{x\} ) }
(1 + |e|)^{\alpha} \right)^4
\right]
 < \infty.
\label{1128c}
\eea
%
%
\end{lemm}
{\em Proof.}
Fix $t, m$ and $x$ with $t \geq 1$, $m/\mu_t \in [1/2,3/2]$,
and  $ x \in tB_0$.
Take cones $C_1, C_2, \ldots, C_K$, 
each with  angular radius $\pi/12$ with 
  point at $x$, and with union $\R^d$,
where $K$ is a constant depending only on $d$. 
For $1 \leq i \leq K$, let $C_i^+$ be the cone
of angular radius $\pi/6$ with point at $x$,
concentric to $C_i$.
Let $1 \leq i \leq K$, and let the random variable $R_i$ be defined as follows:
\begin{itemize}
\item
If there exists a point $Z$ of  $\U_{m,t}$ lying
in the cone $C_i$, and carrying a mark $\TT_Z < \TT_x$,
let $R_i$ be the distance from $x$ to the nearest point $Y$ of  $\U_{m,t}$ lying
in the cone $C_i^+$, and carrying a mark $\TT_Y < \TT_x$.
\item
If no such point $Z$  exists, set $R_i$ to be
$\sup_{y \in tB_0 \cap C_i}|y-x|$, the
furthest distance  from $x$ to any element of
$tB_0 \cap C_i$. 
\end{itemize}
By  Lemma \ref{lem0909}, all incoming edges
to $x$ from points of $\U_{m,t}$ in $C_i$
must be from points at a distance at most $R_i$
from $x$.

Let ${\cal Y}_i$ be the set of points of $\U_{m,t}$
lying in the cone $C_i$ at a distance less than
$R_i$ from $x$. By definition of $R_i$, 
necessarily all
points of $\YY_i$ carry a mark greater than $\TT_x$.
Listing the points of $\YY_i$ as $Y_{i,1},\ldots,Y_{i,\nu(i)}$,
in order of increasing mark,
let $M_i$ be the number of 
 points $Y_{i,j}$ of $\YY_i$ such that $Y_{i,j}$ 
lies closer to $x$ than do any of the points
$Y_{i,1},\ldots, Y_{i,j-1}$ (we include $Y_{i,1}$ in this
set of points). By Lemma \ref{lem0909}, each  
 incoming edge at $x$    with an endpoint
in the cone  $C_i$ is from such a point, so that
\bea
  \left( 
\sum_{ e \in {\cal E}^+(x; \U_{m,t} \cup \{x\} ) }
(1 + |e|)^{\alpha} 
\right)^4
 \leq
\left( \sum_{i=1}^K  
(1 + M_i)(1 + R_i)^{\alpha} 
 \right)^4
\nonumber \\
\leq K^4 
 \sum_{i=1}^K  (1 +M_i)^4(1 + R_i)^{4\alpha}.
\label{1128a}
\eea
When $x$ is inserted into ${\cal U}_{m,t}$,
all removed edges are outgoing from points
that are connected to $x$ after it is inserted.
Hence, the removed edges outgoing from points
 in $C_i$ are the edges outgoing from the points
in the set  $\{Y_{i,1},\ldots,Y_{i,\nu(i)}\}$,
which are connected to $x$ after insertion 
of $x$.
Since these points lie in $C_i \cap B_{R_i}(x)$,
by Lemma \ref{lem0909} 
they lie within distance at most $R_i$ of each other
so that the removed edges outgoing from the points in
$C_i$ are of length at most $R_i$, with the possible sole
exception of the outgoing edge from $Y_{i,1}$. 

Also, we assert that  the removed edge 
 from $Y_{i,1}$  
(if there is one)
has length at most $2 \max_{1 \leq \ell \leq K}R_\ell$.
To see this, note that  if $x$ carries a lower mark
than any point of $\U_{m,t}$, then all of $tB_0$
lies within distance 
 $\max_{1 \leq \ell \leq K}R_i$
of $x$; if not, then for some $\ell \leq K$ there is
a point of $\U_{m,t}$ in $C_\ell^+$ at distance $R_\ell$ from $x$ 
carrying a mark which is lower than  $T_x$, 
and hence also lower than $T_{Y_{i,1}}$; the assertion follows.

By the preceding remarks about removed edges, it follows that
\bea
  \left(\sum_{ e \in {\cal E}^-(x; \U_{m,t} \cup \{x\} ) }
(1 + |e|)^{\alpha} \right)^4
 \leq
 \left( K(1+2  \max_{i \leq K} R_i)^\alpha +
\sum_{i=1}^K    
M_i(1 +  R_i)^{\alpha} \right)^4
\nonumber \\
\leq \left( \sum_{i=1}^K (1 + R_i)^\alpha(2^\alpha K +M_i)
\right)^4 
\nonumber \\
\leq K^4 \sum_{i=1}^K (2^\alpha K +M_i)^4  (1 +R_i)^{4\alpha}.
\label{1128d}
\eea
Conditional on the set of points $\YY_i$ (but not
their marks), any of the $\nu(i)!$ possible 
orderings of the marks of points of $\YY$ is equally likely.
Note that ${M_i \choose 4}$ is the number of
collections of four distinct points $Y_{i,j_1}, \ldots, Y_{i,j_4}$
such that each of $Y_{i,j_k}$, $1 \leq k \leq 4$,
lies closer to $x$ than any point of $Y_{i,1},\ldots, Y_{i, j_k-1}$.

Given $\nu(i)=\ell$,
for any  $j_1< j_2 < j_3 < j_4 \leq \ell$,
the  probability that 
each of $Y_{i,j_k}$, $1 \leq k\leq 4$
lies closer to $x$ than any point of $Y_{i,1},\ldots, Y_{i,j_k-1}$
is equal to $j_1^{-1} j_2^{-1} j_3^{-1} j_4^{-1}$. Hence,
\bean
 \E\left[ \left. {M_i \choose 4} \right|
 \nu (i)=\ell \right] =  \sum_{1\leq i_1 < i_2 < i_3 < i_4\leq \ell } 
\frac{1}{i_1 i_2 i_3 i_4} 
\leq  (1 + \log \ell )^4/4!, ~~~~ \ell >0,
\eean
and since trivially 
$$
M_i^4 \leq 256 {M_i \choose 4} {\bf 1}_{\{ M_i \geq 4\}} 
+ 81 {\bf 1}_{\{ M_i < 4\}} \leq 256 {M_i \choose 4} + 81, 
$$
we obtain
\bean
\E [ M_i^4| 
 \nu(i) = \ell  ] \leq 81 +
 11  (1 + \log \ell)^4, ~~~~~ \ell >0.
\eean
Choose $\eps \in (0,1)$ with $4 \alpha + d \eps < d $.
Conditional on $\nu(i)$, the distribution of $M_i$ does
not depend on the value of $R_i$, so that
for some constant $c_1>0$,
\bean
\E [ M_i^4| \nu(i) , R_i] 
\leq 81 +
 11  (1 + \log \nu(i))^4 {\bf 1}_{\{\nu(i)>0\}} \leq c_1 (1+ \nu(i)^\eps).
\eean

Next, we assert that
the conditional distribution of $\nu(i)$, given $R_i=r$,
is stochastically dominated by the Binomial 
${\rm Bi}(m, \pi_d r^d/ (t^d|B_0|))$,
where for $p >1$ we set ${\rm Bi}(n,p) := {\rm Bi}(n,1)$.
To see this, let $\U_{m,t}^-$ denote the set of points
of $\U_{m,t}$ which carry a mark less than $T_x$,
and let $\U_{m,t}^+:= \U_{m,t} \setminus \U_{m,t}^-$.
Let $N:= \card ( \U_{m,t}^-)$. Then 
$$
\nu(i) = \card( \U_{m,t}^+ \cap B_{R_i}(x) \cap C_i ),
$$
while the value of $R_i$ is determined by the configuration
of $\U_{m,t}^-$. The conditional distribution of $\U_{m,t}^+$,
given $\U_{m,t}^-$, is that of $m-N$ points independently
uniformly distributed in $tB_0$ (thus, this conditional distribution
depends on $\U_{m,t}^-$ only through the value of $N$).
Hence, given $R_i=r$ and $N =n$,
the conditional distribution of  $\nu(i) $ is binomial 
$$
{\rm Bi}(m-n, |B_r(x) \cap C_i \cap tB_0|/|tB_0|),
$$
and since all possible values of $m-N $ are at most $ m$, this
 conditional distribution is stochastically dominated
by
${\rm Bi}(m, \pi_d r^d/ (t^d|B_0|))$,
justifying the assertion above.

By the preceding assertion,
 since we assume $m \leq 2 \lambda t^d |B_0|$,
we have
\bean
\E [ M_i^4| R_i    ] \leq
c_1  \E[ 1 + ({\rm Bi}(m, \pi_d R_i^d/(t^d|B_0|) ) )^\eps]
\\
\leq c_1 (1+ (2 \lambda \pi_d R_i^d)^\eps ) 
\eean
by Jensen's inequality.
Hence,
 for suitable $c_2$,
\bea
\E[ (1+M_i)^4(1+ R_i)^{4\alpha}] \leq 
c_2 
\E[\max (1 , R_i^{4 \alpha + d \eps} )].
\label{1128}
\eea
For any $r>0$, we have  $P[R_i >r ]=0$ unless 
there exists  $y \in tB_0 \cap C_i$ with $|y-x|=r$, in which
case by Lemma \ref{lemgeo}, 
$$
|C_i^+ \cap B_r(x) \cap tB_0| \geq |C_{x,y} \cap B_r(x) \cap tB_0|
\geq c_3 r^d,
$$
for some constant $c_3 >0 $ depending only on $B_0$. 
Hence,
by conditioning on the value of $\TT_x$ 
we have
for large enough $s$ that
\bean
 P[ R_i^{4 \alpha + d \eps} >  s]
\leq \int_0^1 
\left( 1 - \frac{  c_3 u s^{d/(4\alpha + d \eps) } }{t^d|B_0|}  \right)^m du
\\
\leq \int_0^1 \exp( - m c_3 u   s^{d/(4\alpha+ d \eps) }  /(t^d|B_0|)) du
\eean
and since we assume $m \geq \lambda t^d|B_0|/2$,
this is bounded by a constant times
$
 s^{-d/(4\alpha + d \eps) }.  
$
Hence, there is a constant $c_4$ such that
$$
\E[\max(1,R_i^{4 \alpha + d \eps})]
=
 \int_1^\infty P [R_i^{4 \alpha + d \eps} > s] ds 
< c_4.
$$
Hence, using (\ref{1128a}) and  (\ref{1128}),
we obtain (\ref{1128b}). The proof of 
  (\ref{1128c}) using (\ref{1128d}) is similar. $\qed$ \\

\noindent
{\em Proof of Theorem \ref{thOL}.}
We assert that 
 there exists an
 almost surely finite
radius of stabilization $R$  satisfying (\ref{gstab1}) and
(\ref{gstab2}). To see this, 
take a finite
collection of cones $C^+_i$ with point at $\0$ and  angular
radius $\pi/6$, with union $\R^d$;   let  $R_i^*$ be the distance
from $\0$ to the nearest Poisson point in $C^+_i$ to $\0$ with a lower
mark than $T_\0$. It is not hard to see that $R_i^*$
is almost surely finite. Then by Lemma \ref{lem0909}, no
point placed in $C^+_i$ at a distance greater than $R_i^*$ from
$\0$ will be connected to $\0$ in the on-line NNG.
Also, by Lemma \ref{lem0909} again,
 any Poisson point in $C^+_i \cap B_{R_i^*}(\0)$ and
carrying a higher mark than $\0$  has
a lower-marked Poisson point  within distance at most $R_i$,
even before addition of a point  at $\0$, so that its nearest
lower-marked neighbour (before insertion of $\0$) lies
 in $B_{2R_i^*}(\0)$. Hence
the set of edges added or removed upon insertion of a point
 at the origin
is unaffected by changes to $\Po_\la$ outside
$B_{2 \max_i R_i^*}(\0)$; in other words,
$2 \max_i R_i^*$
is a radius of stabilization in the sense of (\ref{gstab1}), (\ref{gstab2}).
%
%
 Thus we can apply
Lemma \ref{lem0516} to get the conditions
(\ref{sstab1}) and (\ref{sstab2}) for either $h(\X,A) = V_\psi^G(\X,A)$
or $h(\X,A)= L^G_\phi(\X,A)$.

The case $\alpha=0$ of Lemma \ref{lemOL} gives us the 
 condition (\ref{ubdmom})  for the functional
$h(\X,A)= V^G_\psi(\X,A)$ for any $\psi \in B(\K_1)$.
Also, the uniform bound (\ref{polybd}) 
is obvious for any such $h$, and by
 scale invariance of the on-line NNG,
 the homogeneity condition
(\ref{eqhomog}) with $\gamma=0$ also holds.
Thus Theorems \ref{1Hclt},  \ref{2Hclt} and \ref{Xncoro} are all applicable
in this case.

Turning to the case where $h(\X,A) = L^G_\phi(\X,A)$, with
$\phi$ satisfying the growth bound (\ref{growthbd}),
once again Lemma \ref{lemOL} gives us the condition (\ref{ubdmom}).
Also, (\ref{polybd}) is again obvious in this case, so that
 Theorems \ref{1Hclt} and  \ref{2Hclt}
are  applicable in this case.  By scale invariance the homogeneity condition
(\ref{eqhomog}) holds (with $\gamma=\alpha$) for the case
$\phi(r) = r^\alpha$, so that Theorem \ref{Xncoro} is also 
applicable in this case.  $\qed$ \\

\section{Proof of the general CLT for lattice systems} 
\label{secprfdscrt}
\allco
Assume throughout this section that $X =(X_z,z \in \Z^d)$ 
 is as described in Section \ref{secmainth}.
Assume also that 
$B_0 \in \RR(\R^d)$ satisfying $|B_0|>0$
 is fixed, and that for $i=1,2,\ldots, k$,
 $(H^i_t(A), t \geq 1, A \in \RR(B_0))$ is
a random set function on $B_0$
 as described in Section \ref{secmainth},  satisfying
the stabilization conditions (\ref{0414a}) and (\ref{0414b})
along with the moments condition (\ref{4moments}) 
for some $\gamma >2$.
Assume also that
 $(t_n)_{n \geq 1}$ is an arbitrary
$[1,\infty)$-valued sequence which tends to infinity as $n \to \infty$.

For $y \in \Z^d$, and $i\in \{1,2,\ldots, k\}$,
since the shifted family of i.i.d. variables $\tau_{y}X$
has the same joint distribution as $X$,
by  (\ref{0414a})  there exists a random variable $\Delta^i_{\infty,y}$ 
such that for $A \in \RR(B_0)$ and $z_n \in \Z^d, n \geq 1$ with 
$\liminf_{ n \to \infty} (\tau_{z_n}(t_n A )) =  \R^d$, 
we have
\bea
 H^i_{t_n,z_n}(\tau_{y}X,A) - H^i_{t_n,z_n}((\tau_{y}X)^\0,A)
 \toP \Delta^i_{\infty,y}.
\label{0426}
\eea
 In other words,
 $\Delta^i_{\infty,y}$ 
is defined in just the same manner
as $\Delta^{H^i}_{\infty}$ at (\ref{0414a}) but using the
 shifted 
family of i.i.d. variables $\tau_{y}X$. 

For $y \in \Z^d$,
let
\bea
F^i_y := \E[ \Delta^i_{\infty,y}|\F_y].
\label{0704a}
\eea
By the conditional Jensen inequality, Fatou's lemma, and
the moments condition (\ref{4moments}), 
\bea
\E[ (F_{y}^i)^2 ] = \E[ ( \E[\Delta^i_{\infty,y} |\F_y] )^2 ]  
\leq 
 \E[ (\Delta^i_{\infty,y})^2 ]  
 = \E[ (\Delta^{H^i}_{\infty})^2 ]  
<\infty.
\label{0622}
\eea
\begin{lemm}
\label{lemerg}
Let $A\in \RR(B_0)$, with $|A|>0$, and let $i, j \in \{1,2,\ldots, k\}$.
 Then
\begin{equation}
(t_n^d|A|)^{-1} \sum_{y \in (t_n A) \cap \Z^d} 
F^i_y F^j_y
  \toone 
 \E[F_\0^i F_\0^j] = 
\E[ \E(\Delta^i_\0(\infty)|\F_\0) \E(\Delta^j_\0(\infty)|\F_\0)
 ].
\label{ergodic}
\end{equation}
\end{lemm}
{\em Proof.}
Since $\F_y$ is the $\sigma$-field generated by
$(\tau_yX_z)_{z \preceq \0}$, the definition of $F_y^i F_y^j$ in terms
of $\tau_yX$ is the same as that of $F_\0^i F_\0^j$ in terms of
$X$. Hence the random field $(F_y^iF_y^j, y\in \Z^d)$
is a stationary family of random variables. Also,
each variable $F_y^i $ has  has finite second moment
by (\ref{0622}),
and likewise for $F_y^j$, so that
$|F_y^iF_y^j|$ has finite first moment by the Cauchy-Schwarz inequality.

Also, the $\sigma$-field of translation-invariant  
 $\sigma(X)$-measurable
events is trivial (see
 Durrett \cite{Dur2}, chapter 6, lemma 4.3).

The result follows from the classical
Ergodic Theorem (\cite{Dur2}, chapter 6, section 2).
For details, see the proof of eqn (2.8) of \cite{Pclt}. 
In the terminology of \cite{Pclt}, the sequence of sets $(t_nA)_{n \geq 1}$
has vanishing relative boundary because of the assumption that
$A$ is Riemann measurable.  This assumption also implies
that $\card(t_nA \cap \Z^d) \sim t_n^d |A|$ as $n \to \infty$.
$\qed$ \\

For $y \in \Z^d$, let  $X^y$
be the random field  $X$ with the value $X_y$
at site $y$ replaced by the independent copy $X_*$ (i.e.,
$X^y = (X^y_z, z \in \Z^d)$ with $X^y_y = X_*$ and $X^y_z
= X_z$ for $z \neq y$).
For $t >0$, and $1 \leq i \leq k$, set  
\bea
\tdelta^i_{t,y}(A) := H^i_{t}(X,A) - H^i_{t}(X^y,A).
\label{dDeldef}
\eea
Observe that $X^y = \tau_{-y}((\tau_yX)^\0)$, so that
\bea
\tdelta^i_{t,y}(A) 
 & = & H^i_{{t},-y}(\tau_{y} X,A) - H^i_{t,-y}((\tau_{y}X)^\0,A).
\label{0424c}
\eea
Therefore by the definition (\ref{deltadef}), since
the translated random field  $\tau_{y}X$ has the same distribution
as $X$,
\bea
\tdelta^i_{t,y}(A) 
& \eqd  & \Delta^{H^i}_{t,-y}(A).  
\label{0424b}
\eea

\begin{lemm}
\label{lem0423}
Let $A \in \RR(B_0)$, and $i \in \{1,2,\ldots, k\}$.
Suppose
$(y_n)_{n \geq 1}$ is  a
$\Z^d$-valued sequence.
 Then 
\bea
\limn
\E [( \tdelta^i_{t_n,y_n}(A) - \Delta^i_{\infty,y_n})^2] =0
~~~~{\rm if}~~
\liminf (\tau_{-y_n}(t_n A ))= \R^d 
\label{0423a}
\eea
and
\bea
\limn
\E [( \tdelta^i_{t_n,y_n}(A)  )^2] =0
~~~~{\rm if}~~
\liminf (\tau_{-y_n}(t_n (B_0 \setminus A )))= \R^d. 
\label{0423b}
\eea
\end{lemm}
{\em Proof.}
The second limiting expression (\ref{0423b}) follows from
the distributional identity (\ref{0424b}) along with
the second stabilization condition (\ref{0414b}) and the moments
condition  (\ref{4moments}) (see \cite{W} A 13.2(f)).

To prove (\ref{0423a}), observe that
since $\Delta^i_{\infty,y}$ is defined in terms of
$\tau_{y}X$ in the same  manner as $\Delta^{H^i}_\infty$
is defined in terms of $X$, 
we have by (\ref{0424c}) that
\bean
\E [( \tdelta^i_{t_n,y_n}(A) - \Delta^i_{\infty,y_n} )^2]
\\
= 
 \E [ ( H^i_{{t_n},-y_n}(\tau_{y_n} X,A) - H^i_{t_n,-y_n}((\tau_{y_n}X)^\0,A) 
- \Delta^i_{\infty,y_n} )^2] 
\\
= 
 \E [ ( H^i_{{t_n},-y_n}(X,A) - H^i_{t_n,-y_n}(X^\0,A) 
- \Delta^{H^i}_{\infty} )^2]. 
\eean
If $\liminf(\tau_{-y_n}(t_n A) )= \R^d$, then this
tends to zero as $n \to \infty$ by 
the stabilization 
and moments 
conditions (\ref{0414a}) 
and (\ref{4moments}) (again see \cite{W} A 13.2(f)).
$\qed$ \\

Recalling the definition of $\widetilde{A}$ at (\ref{1121}), 
define the sequence of sets $(B_n)_{n \geq 1}$ in $\Z^d$ by
\bea
B_n:= \widetilde{t_n B_0} .
\label{Bndef}
\eea
For $i \in \{1,2,\ldots, k\}$, 
 $y \in \Z^d$, $ t \geq 1$, and $A \in \RR(B_0)$, let
\bea
F_{t,y}^i(A) : = \E[ \tdelta^i_{t,y}(A)|\F_y].
\label{0615}
\eea
For $\gamma >1$, the conditional Jensen inequality
implies that
\bea
\E [|F^i_{t,y}(A) |^\gamma] = \E [ |\E [\tdelta^i_{t,y}(A)|\F_y]|^\gamma]
 \leq \E [ |\tdelta^i_{t,y}(A)|^\gamma] 
\label{0426a2} 
\eea
and therefore the distributional identity (\ref{0424b}) 
together with the moments condition (\ref{4moments})
imply that for some $\gamma >2$,
 \bea
\sup\{ 
\E [|F^i_{t,y}(A)|^\gamma] : A \in \RR(B_0), t \geq 1, y \in \widetilde{tB_0}
\} < \infty.
\label{0704}
\eea
\begin{lemm} 
\label{lem0425}
For any $A \in \RR(B_0)$, $A' \in \RR(B_0)$, and any
$i,j \in \{1,\ldots, k\}$,
as $n \to \infty$ we have
\bea
t_n^{-d} \sum_{x \in (t_n(A \cap A') \cap \Z^d) } 
(F^i_{t_n,x}(A) F^j_{t_n,x}(A') - F^i_x F^j_x )  \tolone 0 
\label{0424f}
\eea
and 
\bea
t_n^{-d} \sum_{x \in B_n \setminus t_n(A \cap A') } 
F^i_{t_n,x}(A)F^j_{t_n,x}(A')   \tolone 0 .
\label{0425a}
\eea
\end{lemm}
{\em Proof.}
By the triangle and Cauchy-Schwarz inequalities,
\bea
\E[ | F^i_{t,x}(A)F_{t,x}^j(A') - F^i_xF^j_x | ] 
\leq 
(\E [F_{t,x}^i(A)^2])^{1/2} 
(\E [(F_{t,x}^j(A') - F_x^j)^2])^{1/2} 
\nonumber
\\
+ (\E [ (F^i_{t,x}(A) - F^i_x)^2])^{1/2}
(\E [( F_x^j)^2])^{1/2}. 
~~~~~~~~~~~~~~~~~~  \label{0724c}
\eea

By (\ref{0704}) and (\ref{0622}), $ \E[ (F_{t,x}^i(A))^2 ] $
 and $\E[ (F_{x}^j)^2 ] $ are uniformly bounded.
Moreover, by definitions (\ref{0704a}), (\ref{0615}) 
 and by the conditional Jensen inequality,
\bea
\E[ (F_{t,x}^i(A)-F_x^i)^2 ] = \E[ ( \E[\tdelta^i_{t,x}(A) -
 \Delta^i_{\infty,x}|\F_x] )^2 ]  
\nonumber \\
\leq
 \E[  \E[(\tdelta^i_{t,x}(A) - \Delta^i_{\infty,x})^2|\F_x] ] =
 \E[ (\tdelta^i_{t,x}(A) - \Delta^i_{\infty,x})^2 ]  
\label{0424a}
\eea
and similarly,
\bea
\E[ (F_{t,x}^j(A')-F_x^j)^2 ]
\leq 
 \E[ (\tdelta^j_{t,x}(A')-\Delta_{\infty,x}^j)^2 ].  
\label{0622a}
\eea

For $A \in \RR(B_0)$, define `interior' and `exterior' lattice
sets for the set $t_n A$ by
$$
{\rm int}_n(A) = \{z \in \Z^d: B_{t_n^{1/2}}(z) \subseteq t_n A\};
~~~
{\rm ext}_n(A) = {\rm int}_n(B_0 \setminus A),
$$
and the `boundary' lattice set 
$$
\partial_n(A) =  B_n \setminus ({\rm int}_n (A) \cup {\rm ext}_n(A)),
$$
which consists of lattice points near the boundary either of $t_nA$ or
of $t_n(B_0 \setminus A)$. 

We assert that
\begin{equation}
\limn \sup_{x \in {\rm int}_n(A \cap A')}
 \E[ | 
F_{t_n,x}^i(A)
F_{t_n,x}^j(A')
 - F_x^i F_x^j | ] =0.
\label{limsupeq}
\end{equation}
Indeed, if this were untrue we could take a sequence
$(x_n)_{n \geq 1}$ with $x_n \in {\rm int}_n (A \cap A')$ and 
$$
\limsup  \E[ | F^i_{t_n,x_n}(A)F^j_{t_n,x_n}(A') - 
F^i_{x_n}F^j_{x_n} | ] >0.
$$
 This would imply by (\ref{0724c}), (\ref{0424a}) and (\ref{0622a}) 
that
$$
\limsup_{n \to \infty}  
\max ( \E[  (\tdelta^i_{t_n,x_n}(A) - \Delta^i_{\infty,x_n})^2  ] ,
 \E[  (\tdelta^j_{t_n,x_n}(A') - \Delta^j_{\infty,x_n})^2  ] ) 
>0,
$$
which contradicts eqn (\ref{0423a}) from  Lemma \ref{lem0423}.

By a similar argument to the proof of (\ref{limsupeq}), this time using
(\ref{0423b}) and (\ref{0426a2}), we obtain 
\bea
\limn \sup_{x \in {\rm ext}_n( A)}
 \E[  F_{t_n,x}^i(A)^2  ] =0, ~~~
\limn \sup_{x \in {\rm ext}_n( A')}
 \E[  F_{t_n,x}^j(A')^2  ] =0.
\label{0424e}
\eea
By (\ref{0424e})
 and the Cauchy-Schwarz inequality,
\bea
\lim_{n \to \infty} \sup_{x \in {\rm ext}_n (  A) 
\cup {\rm ext}_n (A')  } 
\E[ |F_{t_n,x}^i( A)  F_{t_n,x}^j(A') |]  =0.
\label{0425b}
\eea
Using the uniform boundedness of
both $ \E[ | F_{t_n,x}^i(A)F_{t_n,x}^j(A') - F_x^i F_x^j | ] $ 
and 
\linebreak 
$ \E[ | F_{t_n,x}^i(A)F_{t_n,x}^j(A') | ] $ (see (\ref{0704}) and
(\ref{0622}))
we may deduce (\ref{0424f}) 
from (\ref{limsupeq}),
and (\ref{0425a}) from (\ref{0425b}).
We here elaborate only on the argument for (\ref{0425a}). 
The  absolute value of the
sum in the left hand side of  (\ref{0425a})
is bounded by four terms, namely
a sum over
 $x \in {\rm ext}_n(A)$, a sum over 
 $x \in {\rm ext}_n(A')$, 
a sum over $x$ in a subset of $\partial_n(A)$, and a sum
over $x$ in a subset of $\partial_n(A')$. The first
two of these terms tend to zero by (\ref{0425b}), while
the other terms tends to zero by the uniform 
boundedness of the terms in the sum and the fact that
the number of sites in $\partial_n(A)$ is small relative to $t_n^d$
(by Riemann measurability of $A$ and of $B_0 \setminus A$),
 and likewise for $A'$.
$\qed$ \\

{\em Proof of Theorem \ref{basictheo}.}
We consider linear combinations.  Recall that $(t_n)_{n \geq 1}$
is an arbitrary sequence tending to infinity, and
let $b_1,\ldots, b_k$ be arbitrary constants.
By the Cram\'er-Wold device (see, e.g., \cite{Dur2}) 
it suffices to prove that with  
\bea
\sigma^*_{j\ell} := 
\E( \E[\Delta^{H^j}_\infty|\F_\0] \E[\Delta^{H^\ell}_\infty|\F_\0]),
\label{0724d}
\eea
 we have 
\bea
 t_n^{-d/2} \sum_{j=1}^k  b_j( H^j_{t_n}(A_j) - \E H^j_{t_n}(A_j) ) \tod 
\NN \left(0, \sum_{j=1}^k\sum_{\ell=1}^k b_jb_\ell |A_j \cap A_\ell| 
\sigma^*_{j\ell} \right),
\label{0424d}
\eea
and that  the variance of the left hand side of (\ref{0424d})
 converges to that of the right hand side.
  We shall represent the left hand side of (\ref{0424d}) 
as a sum of martingale differences.

For $j,\ell \in \{1,2,\ldots, k\}$,
let $A_{n,j,\ell} := t_n( A_j \cap A_\ell) \cap \Z^d$.
Recall that $B_n := \widetilde{t_n B_0}$ (see (\ref{Bndef})). 
 Let $\nu_n = \card(B_n)$ and
$\nu_{n,j,\ell}:= \card(A_{n,j,\ell})$.
Since $B_0,A_1,A_2,\ldots, A_k$ are all 
Riemann measurable we have (for each $j, \ell$)
\bea
\lim_{n\to \infty} ( \nu_n/t_n^d) = |B_0|; ~~~~~
\lim_{n\to \infty} ( \nu_{n,j,\ell}/t_n^d) = |A_j \cap A_\ell|.
\label{0424}
\eea
Define the filtration $(\G_0,\G_1,\ldots,\G_{\nu_n})$  as follows: 
 let $\G_0$ be the trivial $\sigma$-field, 
label the elements of $B_n$ in lexicographic order as
$x_1,\ldots,x_{\nu_n}$, and let
$\G_i = \F_{x_i}$ for $1 \leq i \leq \nu_n$.
Then
$$
\sum_{j=1}^k b_j (H^j_{t_n}(A_j) - \E H^j_{t_n}(A_j) ) 
= \sum_{i=1}^{\nu_n} D_i, 
$$
where we set $D_i := \sum_{j=1}^k b_j D_{i,j} $ with
$$
D_{i,j}: =
 \E[H^j_{t_n}(A_j)|\G_i] -  \E[H^j_{t_n}(A_j)|\G_{i-1}] .
$$
By orthogonality of martingale differences, 
$$
\Var \left[ \sum_{j=1}^k b_j H^j_{t_n}(A_j) \right]
  = \E \sum_{i=1}^{\nu_n} D_i^2.
$$
By this representation of the variance, along with 
 the central limit theorem for martingale difference arrays
 (Theorem (2.3) of McLeish \cite{McL}, or Theorem  2.10 of Penrose \cite{Pbk})
 it suffices to prove the conditions
\begin{equation}
\label{1MCT1}
\sup_{n \geq 1}
\E \left[ \max_{1 \leq i \leq \nu_n} \left(  t_n^{-d/2} |D_i| \right)^2 \right]
< \infty , 
\end{equation}
\begin{equation}
t_n^{-d/2} \max_{1 \leq i \leq \nu_n} |D_i| \toP 0,
\label{1MCT2}
\end{equation}
and 
\begin{equation}
t_n^{-d} \sum_{i=1}^{\nu_n} D_i^2 \toone \sum_{j=1}^k \sum_{\ell=1}^k 
b_jb_\ell |A_j \cap A_\ell| \sigma^*_{j\ell} .
\label{1MCT3}
\end{equation}

With $\tdelta^j_{t,y}(A)$ defined at (\ref{dDeldef}), and
$F^j_{t,y}(A)$ defined at (\ref{0615}), we have
\begin{equation}
D_{i,j} 
 = E [ \tdelta^j_{t_n,x_i}(A_j) |\F_{x_i}] = F^j_{t_n,x_i} (A_j). 
\label{altdi}
\end{equation}
First we check  (\ref{1MCT1}).
By (\ref{altdi}),
we have
$$
t_n^{-d} \E [\max_{i \leq \nu_n} D_i^2 ] \leq
 t_n^{-d} \sum_{i=1}^{\nu_n} \E[D_i^2] = 
 t_n^{-d} \sum_{i=1}^{\nu_n} \E\left[\left(
\sum_{j=1}^k b_j 
F^j_{t_n,x_i}(A_j)
\right)^2 \right]  
$$
which is bounded, uniformly in $n$, by (\ref{0424}) and (\ref{0704}).

For the second condition (\ref{1MCT2}), let $\eps >0$
and use Boole's and Markov's 
inequalities to obtain
\bean
P \left[  \max_{1 \leq i \leq \nu_n} |D_i| \geq t_n^{d/2} \eps \right] 
\leq
\sum_{i=1}^{\nu_n} \frac{\E[|D_i|^\gamma]}{ t_n^{\gamma d/2} \eps^\gamma },
\eean
which tends to zero,  by 
 (\ref{0424}) and the fact that for some $\gamma >2$,
 $\E [|D_i|^\gamma]$
is bounded, uniformly over $n \geq 1$ and $i \leq \nu_n$,
 by (\ref{altdi}) and (\ref{0704}).

It remains to prove (\ref{1MCT3}). 
It suffices to prove that for each $j, \ell \in \{1,2,\ldots,k\}$ 
we have
\bea
t_n^{-d} \sum_{i=1}^{\nu_n} D_{i,j}D_{i,\ell}
 \toone |A_j\cap A_\ell| \sigma^*_{j \ell} . 
\label{030415}
\eea

Using  (\ref{ergodic}), (\ref{0424f}), (\ref{0425a}),
and (\ref{0424}), we obtain
\bea
t_n^{-d} 
\sum_{x \in \tnlB } F^j_{t_n,x}(A_j)F^\ell_{t_n,x}( A_\ell)  \toone 
 |A_j \cap A_\ell| \E[F_\0^j F_\0^\ell]. 
\label{0724e}
\eea
By the definitions (\ref{0704a}) and (\ref{0724d}),
the right-hand side of (\ref{0724e}) equals 
$|A_j \cap A_{\ell}| \sigma^*_{j\ell}$.
By (\ref{altdi}),
 eqn (\ref{0724e}) 
 gives us (\ref{030415}).
The proof is complete.
$\qed$

\section{Proof of general continuum results}
\label{secprfcont}
\allco
In this section we prove 
the results stated in Section \ref{contsubsec}.
Recall the definition of a point process set function at the start of
Section \ref{contsubsec} and the definition of add one cost
$\delta(A,\X)$ given at (\ref{add1c}). 
First we give some consequences of the stabilization and moments
conditions given in that section. 

Given $\la >0$, given a point process set function $\HH$ and
a random variable $\delinf$, 
let us say $h$ is {\em weakly stabilizing} at intensity $\la$ 
with stabilizing limit $\delinf$, if  for any  $B_0 \in \RR(\R^d)$
with $|B_0|>0$, for any $A \in \RR(B_0)$, 
and  any
 $([1,\infty) \times \R^d)$-valued sequence $(t_n,x_n)_{ n \geq 1}$, we have
\bea
\delta (
\tau_{x_n}(t_nA),\Po_\la \cap (\tau_{x_n} (t_nB_0))
 ) \toas \delinf ~~~~{\rm 
if} ~~ \liminf_{ n \to \infty} (\tau_{x_n}(t_n A ) )=  \R^d
\label{0513a}
\eea
and
\bea
\delta ( \tau_{x_n}(t_nA),\Po_\la \cap (\tau_{x_n} (t_nB_0))) 
\toas  0 ~~~~{\rm 
if} ~~ \liminf_{ n \to \infty} (\tau_{x_n}(t_n (B_0 \setminus A) ) )=  \R^d.
\label{0513b}
\eea
Weak stabilization can be viewed as a
 continuum version of the conditions
(\ref{0414a}), (\ref{0414b}). We also consider 
$h$ satisfying the moments condition
\bea
\sup_{A \in \RR(B_0), t \geq 1, x \in \R^d : \0 \in \tau_x(tB_0)}  
\{ \E [ \delta (\tau_x(tA),\Po_\la \cap \tau_x(tB_0))^4] \}
 < \infty .
\label{pbdmom}
\eea
which is a Poisson point process version of the condition 
(\ref{4moments}).

\begin{lemm}
\label{lemstrongweak}
Let $\la >0$. Suppose 
$\HH$ is a  point process set function. Then:

(i) If $\HH$ is strongly stabilizing 
at intensity $\la$ with stabilizing limit $\delinf$
(i.e., satisfies  
$(\ref{sstab1})$ and $(\ref{sstab2}))$, then  
$\HH$ is weakly stabilizing 
at intensity $\la$ with stabilizing limit $\delinf$
(i.e., satisfies $(\ref{0513a})$ and $(\ref{0513b}))$.


(ii) If $h$ satisfies conditions $(\ref{ubdmom})$ and $(\ref{polybd})$,
then $h $  satisfies $(\ref{pbdmom})$.
\end{lemm}
{\em Proof.}
Part (i) is obvious.  Part  (ii) is proved by
 a similar argument to the proof of Lemma 4.1 of \cite{PY}, which
we omit. $\qed$ \\

Suppose that $\la >0$ and 
we are given a 
point process set function $\HH$ that
is weakly stabilizing at intensity $\la$ with stabilizing limit
$\delinf$. 
For any locally finite set $\X \subset \R^d$ and
any $x \in \R^d$, define
\bea
\delta_\infty(x,\X) : = \limsup_{n \to \infty} (
h[ (\X \cap B_n(x)) \cup \{x\}, B_n(x) ] - h[ \X \cap B_n(x),B_n(x)]). 
\label{def0911}
\eea
By the definition (\ref{add1c}) of add one cost, and
 translation invariance (\ref{transinv}), we have
$$
\delta_\infty(x,\Po_\la) : = \limsup_{n \to \infty}
 \delta(B_n(\0) , (\tau_{-x}(\Po_\la)) \cap B_n(\0) ). 
$$
Since $\tau_{-x}(\Po_\la)$ is a homogeneous
Poisson process of intensity $\la$,
by taking $A = B_0 = B_1(\0)$, $x_n =\0$ and $t_n =n$ in (\ref{0513a})
 we see that
 $\delta_\infty(x,\Po_\la)$   almost surely
 equals the stabilizing limit
$\delta_\infty(\tau_{-x}(\Po_\la))$ of $h$ with respect to
the shifted Poisson process $\tau_{-x}(\Po_\la)$.
Thus for all $x \in \R^d $, we have as $n \to \infty$ that
\bea
h( (\Po_\la \cap B_n(x)) \cup \{x\}, B_n(x) ) - 
h( \Po_\la \cap B_n(x),B_n(x))   
\toas \delta_\infty(x,\Po_\la).
\label{from0911}
\eea

\begin{lemm}
\label{lem0912}
Let $\la >0$ and let the point process set function
$h$ be weakly stabilizing
at intensity $\la$.
Given $\eps \in (0,1]$, let
 the random vector $\xi(\eps)$ be uniformly distributed
over the cube $Q_\0^\eps$, independent of $\Po_\la$. Then
\bea
\delta_\infty(\xi(\eps),\Po_\la)
 \toP \delinf ~~~{\rm as} ~~\eps \downarrow 0.
\label{0907}
\eea
\end{lemm}
{\em Proof.}
As $K \to \infty$, we have by (\ref{0513a}) and 
(\ref{from0911}) that
\bean
h((\Po_\la \cap B_K(\0)) \cup \{\0\}, B_K(\0)) - h(\Po_\la  \cap
B_K(\0), B_K(\0)  ) \toas \delinf ;  
\\
 h((\Po_\la  \cap B_K(\xi(\eps))) \cup \{\xi(\eps)\}, B_K(\xi(\eps)))
~~~~~~~~~~~~~
~~~~~~~~~~~~~
\\
 - h(\Po_\la  \cap B_K(\xi(\eps)), B_K(\xi(\eps)) )
 \toas
\delta_\infty(\xi(\eps),\Po_\la ).
\eean 
Also, for any $K$ it is the case that $P[ \Po_\la \cap B_K(\xi(\eps)) \neq
\Po_\la \cap B_K(\0) ] \to 0$ as $\eps \downarrow 0$.  
Hence, it suffices to prove that for any integer $j \geq 1$, and any $K>0$, 
if $X_1,\ldots X_j $ are
 uniformly distributed over $ B_K(\0)$, 
 independent of each other and of $\xi(\eps)$,
then
\bea
h(\{\xi(\eps),X_1, \ldots, X_j\},B_K(\xi(\eps)))
 \toP 
 h(\{\0,X_1, \ldots, X_j\},B_K(\0)) 
~~~{\rm as} ~~\eps \downarrow 0;
\label{0912a3}
\\
 h(\{X_1, \ldots, X_j\},B_K(\xi(\eps))) \toP  h(\{X_1, \ldots, X_j\},B_K(\0))
~~~{\rm as} ~~\eps \downarrow 0.
\label{0912a2}
\eea
By (\ref{polybd}), the above random variables
are uniformly bounded by a constant (dependent on $j$ and $K$). 
Define
$\tilde{h}: (\R^d)^{j} \to \R$ and
$h^*: (\R^d)^{j} \to \R$ 
by
\bean
\tilde{h} (y_1, \ldots, y_j) := 
h(\{\0,y_1, y_2 + y_1 ,y_3+y_1, \ldots, y_j +y_1\},B_K(\0)) ;
\\
h^* (y_1, \ldots, y_j) := 
h(\{y_1, y_2 + y_1 ,y_3+y_1, \ldots, y_j +y_1\},B_K(\0)) . 
 \eean
If $x_1$ lies at a Lebesgue point (see e.g. \cite{Rud2}) of
$\tilde{h}(\cdot, x_2-x_1, x_3-x_1, \ldots, x_j-x_1)$
then
\bea
\eps^{-d} \int_{Q_\0^\eps} 
|h(\{x,x_1, x_2,\ldots,x_j\},B_K(x)) - h(\{\0,x_1,x_2,\ldots,x_j\},B_K(\0))| dx
\nonumber \\
= \eps^{-d} \int_{Q_\0^\eps} 
|h(\{\0,x_1-x, x_2-x,\ldots,x_j-x\},B_K(\0)) ~~~~~~~ 
\nonumber \\
- h(\{\0,x_1,x_2,\ldots,x_j\},B_K(\0))| dx
\nonumber \\
= \eps^{-d} \int_{Q_\0^\eps} 
|\tilde{h}(x_1-x, x_2-x_1,\ldots,x_j-x_1) - \tilde{h}
(x_1,x_2-x_1,\ldots,x_j-x_1)| dx
\nonumber \\
\to 0 {\rm ~~ as ~~} \eps \downarrow 0,
~~~~~~~~~~~~~~
\label{0912}
\eea
where the last line comes from the  definition
of a Lebesgue point. Similarly, if 
 $x_1$ lies at a Lebesgue point  of
$h^*(\cdot, x_2-x_1, x_3-x_1, \ldots, x_j-x_1)$, then
\bea
\eps^{-d} \int_{Q_\0^\eps} 
|h(\{x_1, x_2,\ldots,x_j\},B_K(x)) - h(\{x_1,x_2,\ldots,x_j\},B_K(\0))| dx
\nonumber \\
\to 0 {\rm ~~ as ~~} \eps \downarrow 0.
\label{0916}
\eea

Since we assume $h(\{x_1,\ldots,x_j\},B_K(\0))$ is a Borel-measurable function
of $(x_1,\ldots,x_j)$, it follows that 
for all $(y_2, \ldots,y_j) \in (\R^d)^{j-1}$ the function
$\tilde{h}(\{\cdot,y_2, \ldots,y_j\})$ is Borel-measurable, and
hence, by the Lebesgue Density Theorem (see \cite{Rud} or
\cite{Pbk}),
that almost every $x \in \R^d$ is a Lebesgue point of
$\tilde{h}(\cdot, y_2, y_3, \ldots, y_j)$.

Suppose $X_1, \ldots, X_j$ are independent and
 uniformly distributed over $B_K(\0)$. Then 
for almost every possible collection of values for
 $(X_2-X_1,\ldots,X_j-X_1)$ 
the conditional distribution  of $X_1$ 
conditional on these values of $(X_2-X_1,\ldots,X_j-X_1)$ 
is absolutely continuous with respect
to Lebesgue measure on $\R^d$ (in fact, uniform over
a certain region).  Hence, given the values of
$(X_2 -X_1 , \ldots,X_j -X_1)$, the conditional probability
that $X_1$ lies at a Lebesgue point of 
$\tilde{h}(\cdot, X_2-X_1, X_3-X_1, \ldots, X_j-X_1)$ is 1.
Thus, with probability 1, 
 $X_1$ lies at a Lebesgue point of 
$\tilde{h}(\cdot, X_2-X_1, X_3-X_1, \ldots, X_j-X_1)$.
Hence by (\ref{0912}), and the Dominated Convergence Theorem,
\bean
\E[ |h(\{\xi(\eps),X_1,\ldots, X_j\},B_K(\xi(\eps))) - h(\{\0,X_1,\ldots, X_j\},B_K(\0)) | ]
\\
=  
\E \left[ \eps^{-d} \int_{Q_\0^\eps}
 |h(\{x,X_1,\ldots, X_j\},B_K(\xi(\eps))) - h(\{\0,X_1,\ldots, X_j\},B_K(\0)) |
 dx \right] 
\\
\to 0 {\rm ~~as~~ \eps \downarrow 0}.
\eean
Since convergence in $L^1$ implies convergence in probability,
  (\ref{0912a3}) then  follows. Also, by a similar argument to
the above, $X_1$ lies almost surely at a Lebesgue
point  of $h^*(\cdot,X_2-X_1,X_3-X_1,\ldots, X_j-X_1)$
so that using (\ref{0916}) we obtain convergence in
$L^1$ of $h(\{X_1,\ldots, X_j\},B_K(\xi(\eps))) $ to
$ h(\{X_1,\ldots, X_j\},B_K(\0)) $, to obtain (\ref{0912a2}).
$\qed$\\

Next we use
 discretization and application of Theorem \ref{basictheo} to prove
a weaker statement of
Theorem \ref{1Hclt}, which does not
include the expression (\ref{0905}) for $\sigma_{ij}^\la$.
In the proof we introduce a parameter $\eps$ which
we shall later on make tend to zero to establish (\ref{0905}).

\begin{prop}
\label{prop1121}
Let $\la >0$ and let 
$B_0 \in \RR(\R^d)$ with $|B_0|>0$.
  Suppose that $\HH^1, \ldots, \HH^k$ 
are point process set functions which satisfy the weak
 stabilization conditions
 $(\ref{0513a})$, $(\ref{0513b})$, 
and  the moments condition 
 $(\ref{pbdmom})$. 
 Then there exists a $k \times k$  matrix
 $\Sigma^\la =(\sigma_{ij}^\la)_{i,j=1}^k $ such that
if $A_1,\ldots, A_k$ are sets in $\RR(B_0)$, then
as $t \to \infty$,  
\bea
t^{-d} \cov(h^i(\Po_{\la,t},tA_i), h^j(\Po_{\la,t},tA_j))
\to \lambda \sigma_{ij}^\la|A_i \cap A_j|
\label{1121d}
\eea
 and
\bea
t^{-d/2} (\HH^i(\Po_{\la,t},A_i) - \E \HH^i( \Po_{\la,t},A_i) )_{i=1}^k \tod 
\NN(\0,  (\lambda \sigma_{ij}^\la|A_i \cap A_j|)_{i,j=1}^k).
\label{1121e}
\eea
\end{prop}
{\em Proof.}
Fix $\eps \in (0,1]$.
To apply Theorem \ref{basictheo}, for $z \in \Z^d$ define $X_z$ to be the
point process $\tau_{-\eps z}(\Po_\la \cap Q_z^\eps)$.  Then $X_z$
 ($z \in \Z^d$) are independent and identically distributed
(they are independent Poisson processes on $Q_\0^\eps$ of intensity
$\la$). Also, define
the random set function
\bea
H^i_t(A) :=  h^i( \Po_{\la,t} , t\eps A), ~~~~ t \geq 1,  A \in \RR(\eps^{-1}B_0),  
\label{0910}
\eea
which is a function of $(X_z, z \in \widetilde{(t/\eps)B_0})$; here
 we denote this function by $g((X_z,   z \in \widetilde{(t/\eps)B_0}))$.
Set $H^i_{t,y}:= H^i_t(\tau_yX,A)$, as at (\ref{Htzdef}).
Then $H_{t,y}^i = g((X_{y+z}, z \in \widetilde{(t/\eps)B_0}))$,
and hence 
by the translation invariance property (\ref{transinv})  of $h^i$ we have
\bea
H^i_{t,y}(A) =
 h^i (\tau_{-\eps y}( \Po_{\la}) \cap  tB_0  ,t\eps A)
 = h^i (  \Po_\la \cap \tau_{\eps y}(tB_0),\tau_{\eps y} (t\eps A)).
\label{0622b}
\eea
[For example, if $h^i(\X;A)$ is simply the number of points of $\X$
in $A$, then (using the definition of $X_z$ above) 
$
g((X_z,   z \in \widetilde{(t/\eps)B_0})) $ equals $ \sum_z \card( \tau_{\eps z}(X_z)
\cap t \eps A ), 
$
and hence, 
\bean
 g((X_{y+z}, z \in \widetilde{(t/\eps)B_0})) = \sum_z \card
( \tau_{\eps z}(X_{y+z}) \cap t \eps A) 
\\
= \sum_z \card ( \tau_{-\eps y}(\Po_\la \cap Q^\eps_{z+y}) \cap t \eps A )
= \card ( \tau_{-\eps y}(\Po_\la) \cap t \eps A ),
\eean  
which is consistent with (\ref{0622b}).]

We need to check conditions (\ref{0414a}), (\ref{0414b}), 
and (\ref{4moments}) in this context. These refer to the increment
$$
\Delta_{t,y}^{H^i}(A) = H^i_{t,y}(X,A) - H^i_{t,y}(X^\0,A)
$$
which, by (\ref{0622b}), is (minus) the increment in
$h^i(\Po_\la \cap \tau_{\eps y}(tB_0),\tau_{\eps y}(t\eps A))$
when we resample the  Poisson process $\Po_\la$ in the cube
$Q_\0^\eps$. The stabilization condition (\ref{0513a}) refers instead
to the insertion of a single point at the origin; however, the
required stabilization (\ref{0414a}) (in the present Poissonian
context) can be deduced from (\ref{0513a}) by the argument used to
prove (3.2) of \cite{PY}. Moreover, a virtually
identical argument can be used to deduce 
 (\ref{0414b}) from (\ref{0513b}).

The proof of  (\ref{4moments}), in this context, from
the assumed condition (\ref{pbdmom}), proceeds 
essentially by the argument given to prove (3.3) of
\cite{PY}; because of this proximity we do not give further details.
Having established conditions (\ref{0414a}), (\ref{0414b}) and (\ref{4moments})
we may  apply Theorem \ref{basictheo} to deduce the results
(\ref{1121d}) and (\ref{1121e}) (see (\ref{0910a})
 and (\ref{1910b}) below).  $\qed$ \\

The proof of proposition \ref{prop1121} just given actually provides
us with some information about the limiting variance matrix
$(\sigma_{ij}^\la)_{i,j=1}^k$.
In the context of  this proof,
the $\sigma$-field $\F_\0$ appearing in
 Theorem \ref{basictheo} is, in effect,
 the $\sigma$-field 
 generated by
the restriction of the Poisson configuration $\Po_\la$
to $\cup_{z \in \Z^d, z \preceq \0} Q_z^\eps$, i.e.,
to cubes in the division of $\R^d$ into cubes $Q_z^\eps$
of side $\eps$, up to and including $Q_\0^\eps$ in the lexicographic
ordering.
To emphasize its dependence on $\eps$, 
 we denote this $\sigma$-field
by $\F_\0^\eps$.
With the random set function $H_t^i$ defined by 
(\ref{0910}), define $\sigma^*_{ij}(\eps)$ (which also depends on $\la$)
by
\bea
\sigma_{ij}^*(\eps) := 
 \E[ \E(\Delta^{H^i}_\infty|\F_\0^\eps)
 \E(\Delta^{H^j}_\infty|\F_\0^\eps)  ]. 
\label{0906}
\eea
Then the application of Theorem \ref{basictheo} in the preceding proof gives us 
\bea
\lim_{t \to \infty} t^{-d} \cov ( 
h^i(\Po_{\la,t}, tA_i), h^j(\Po_{\la,t}, tA_j) ) =
\lim_{t \to \infty} t^{-d} \cov ( 
H_t^i(\eps^{-1} A_i), H_t^j(\eps^{-1} A_j) )  )
\nonumber
\\
= \eps^{-d} |A_i \cap A_j| \sigma_{ij}^* (\eps) 
~~~~~~~~~~
\label{0910a}
\eea 
and as $t \to \infty$,
\bea
(t^{-d/2} ( h_t^i(\Po_{\la,t},tA_i) - \E h_t^i(\Po_{\la,t},tA_i) 
))_{i=1}^k \tod 
\NN(\0, (\eps^{-d}\sigma^*_{ij}(\eps) |A_i \cap A_j|)_{i,j=1}^k ).
\label{1910b}
\eea
In other words,
the matrix $(\sigma_{i,j}^\la)_{i,j=1}^k$
in the statement of  Proposition \ref{prop1121} is given,
 for any $\eps \in (0,1]$,
by
\bea
\sigma_{ij}^\la = \lambda^{-1} \eps^{-d} \sigma^*_{ij}(\eps). 
\label{0910b}
\eea
{\em Proof of Theorem $\ref{1Hclt}$.}
In view of Lemma \ref{lemstrongweak},
 Proposition \ref{prop1121} and the discussion above,
it remains to prove that if $h^1,\ldots, h^k$ are
weakly stabilizing at intensity $\la$, and satisfy the moments condition
(\ref{pbdmom}), then
 $\sigma_{ij}^\la$, given
by (\ref{0910b}) for any $\eps \in (0,1]$,  is also given by (\ref{0905}).

With $i$ and $j$ fixed, define 
point process set functions $h:= h^i +h^j$ and
$h':= h^i - h^j$, along with 
the corresponding random set functions
$H := H^i + H^j$ and $H' := H^i - H^j$ (where $H^i$ and $H^j$ are given at 
(\ref{0910})).
  The definition of $H^i$
also depends on $\eps$, as does the limiting increment $\Delta_\infty^{H^i}$;
from now on we denote the latter quantity
by $\Delta_\infty^{H^i,\eps}$, and define 
$\Delta_\infty^{H,\eps}$ and  $\Delta_\infty^{H',\eps}$  analogously.
By linearity, for all $\eps >0$ we have
\bean
\sigma_{ij}^*(\eps) = (1/4) 
 \E\left[ (\E[\Delta^{H,\eps}_\infty|\F_\0^\eps] )^2
 - (\E[\Delta^{H',\eps}_\infty|\F_\0^\eps ] )^2 \right].
\eean
To prove (\ref{0905}), we use the fact that the value of 
 $\eps^{-d}\sigma^*_{ij}(\eps)$ does not depend on the choice of
 $\eps $, since the left hand side of (\ref{0910a})
does not depend on $\eps$ and therefore neither does the
right hand side. 
The aim is to show, by taking $\eps \downarrow 0$ in (\ref{0910b})
 that $\sigma_{ij}^\la$ equals the expression
\bean
\E[ \E(\delinfi| \F) \E(\delinfj|\F)]
= \frac{1}{4} \: 
\E \left[ (\E[\delinf| \F])^2 - (\E[\delinfd|\F])^2\right],
\eean
where $\delinf$ (respectively $\delinfd$)  is the stabilizing
limit of the point process set function $h^i + h^j$ 
(respectively $h^i - h^j$).
In  other words, it remains to prove that
\bea
\lim_{\eps \downarrow 0} \eps^{-d}
\E[ (\E[\Delta^{H, \eps}_\infty|\F_\0^\eps])^2] 
= \lambda \E[ (\E[\delinf| \F])^2 ],
\label{0906b}
\eea
and also a  similar limit 
for $H'$, for which the proof will be identical.

By following the proof of (\cite{PY}, Lemma 3.1)
and observing that $c(\mu)$ in that proof tends to zero as $\mu \downarrow 0$,
  we see that
\bea
\lim_{\eps \downarrow 0} \E[(\Delta^{H, \eps}_\infty)^4 ] =0.
\label{0906a}
\eea
Let $N_\eps$ (respectively $N'_\eps$) be the number of
points of $\Po_\la$ in the cube $Q_\0^\eps$ 
(respectively the number of resampled Poisson points in $Q_\0^\eps$).
If $N_\eps=N'_\eps=0$ then $\Delta^{H,\eps}_\infty =0$. Also,
$N_\eps$ is $\F_\0^\eps$-measurable. Hence, 
\bean
\E[\Delta_\infty^{H,\eps} {\bf 1}_{\{N'_\eps=0\}} | \F_\0^\eps ]
 {\bf 1}_{\{N_\eps=0\}}
=
\E[\Delta_\infty^{H,\eps} {\bf 1}_{\{N'_\eps=0\}}
{\bf 1}_{\{N_\eps=0\}}
 | \F_\0^\eps ] 
\\
= \E [ 0| \F_\0^\eps] = 0,~~~~{\rm a.s.}
\eean
Hence
\bean
(\E[\Delta_\infty^{H, \eps} | \F_\0^\eps ] {\bf 1}_{\{N_\eps=0\}})^2
=  
(\E[\Delta_\infty^{H, \eps} 
( {\bf 1}_{\{N'_\eps=0\}} + {\bf 1}_{\{N'_\eps>0\}}  ) 
|\F_\0^\eps ] 
{\bf 1}_{\{N_\eps=0\}}
)^2
\\
=
(\E[\Delta_\infty^{H, \eps} {\bf 1}_{\{N'_\eps>0\}}  
 |\F_\0^\eps ] {\bf 1}_{\{N_\eps=0\}})^2
\leq
(\E[\Delta_\infty^{H, \eps} {\bf 1}_{\{N'_\eps>0\}}  |\F_\0^\eps ] )^2,
~~~ {\rm a.s.}.
\eean
  Hence, by the conditional Cauchy-Schwarz
 inequality (see e.g. \cite{Dur2}), and the independence
of $N'_\eps$ and $\F_\0^\eps$,
\bean
(\E[\Delta_\infty^{H, \eps} | \F_\0^\eps ] {\bf 1}_{\{N_\eps=0\}})^2
\leq 
 \E[(\Delta_\infty^{H, \eps})^2 |\F_\0^\eps ]
 P[N'_\eps >0],~~~~{\rm a.s.}  
\eean
Taking expectations, then using Jensen's inequality and
(\ref{0906a}), we obtain
\bea
\eps^{-d} 
\E[
(\E[\Delta_\infty^{H, \eps} | \F_\0^\eps ])^2 {\bf 1}_{\{N_\eps=0\}}]
\leq 
\eps^{-d} 
P[N'_\eps >0]   
\E[(\Delta_\infty^{H, \eps})^2  ]  
\nonumber
\\
 \leq 
\lambda
\E[(\Delta^{H, \eps}_\infty)^4 ]^{1/2}
\nonumber \\
 \to 0 ~~~{\rm as}~ \eps \downarrow 0.  
\eea
Let $Y_\eps= (\E(\Delta^{H, \eps}_\infty|\F_\0^\eps))^2$. 
By the Cauchy-Schwarz and Jensen inequalities, and (\ref{0906a}),
\bea
\eps^{-d} 
\E[Y_\eps {\bf 1}_{\{N_\eps \geq 2\}} ] \leq 
\eps^{-d} (P[ N_\eps \geq 2])^{1/2} (\E[ Y_\eps^2])^{1/2} 
\nonumber \\
\leq {\rm const.} \times 
\E[(\Delta^{H, \eps}_\infty)^4]^{1/2} 
\to 0 {\rm ~~~ as ~~} \eps \downarrow 0.
\eea
Similarly, 
$$
\eps^{-d} \E[Y_\eps {\bf 1}_{\{ N_\eps =1, N'_\eps \geq 1\} } ] 
\to 0 {\rm ~~~ as ~~} \eps \downarrow 0.
$$
It remains to consider $\E[Y_\eps {\bf 1}_{\{ N_\eps=1, N'_\eps=0\}}]$.
Since $P[N_\eps=1, N'_\eps=0] \sim \lambda \eps^d$ as $\eps \downarrow 0$,
to establish (\ref{0906b}) we must show that
\bea
\lim_{\eps \downarrow 0}
\E[ (\E[\Delta_\infty^{H,\eps}|\F_\0^\eps])^2| N_\eps=1, N'_\eps=0 ]
 = \E[(\E[\delinf|\F])^2].
\label{0912d}
\eea
Given $\eps$, let $\Y_\eps$ (respectively $\Y'_\eps$)
 be the restriction of the Poisson 
 process $\Po_\la$  to the union of  cubes $Q_z^\eps$ with
$z \in \Z^d$ and $z \prec \0$  (respectively, with $\0 \prec z$).

Given that $N_\eps=1$ and $N'_\eps=0$, the restriction
of $\Po$ to $Q_\0^\eps$ consists 
of a single point  uniformly distributed over $Q_\0^\eps$ 
and independent of $(\Y_\eps,\Y'_\eps)$;
we denote this random point by $\xi'(\eps)$. Then, 
given that $N_\eps=1$ and $N'_\eps=0$, 
almost surely $\Delta_\infty^{H,\eps} $ equals the increment 
$\delta_\infty(\xi'(\eps),\Y_\eps \cup \Y'_\eps)$
(using notation defined at (\ref{def0911})).  Thus, 
\bea
\E[ (\E[\Delta_\infty^{H,\eps}|\F_\0^\eps])^2| N_\eps=1, N'_\eps=0 ]
= \E [  (\E[\delta_\infty(\xi(\eps), \Y_\eps \cup \Y'_\eps)|\xi(\eps),\Y_\eps])^2] 
\label{0912c}
\eea
where, as in Lemma \ref{lem0912},
 $\xi(\eps)$ is uniformly distributed over $Q_\0^\eps$
and is independent of $\Po_\la$.

By the Cauchy-Schwarz  and Jensen inequalities,
\bea
\E \{ ( \E[\delta_\infty(\xi(\eps),\Y_\eps \cup \Y'_\eps)|\xi(\eps),\Y_\eps])^2
- (\E[ \delinf|\xi(\eps),\Y_\eps])^2 \}
\nonumber \\
= \E  \{  
\E[\delta_\infty(\xi(\eps),\Y_\eps \cup \Y'_\eps) + \delinf|\xi(\eps),\Y_\eps] 
~~~~~~~~
\nonumber \\
\times
\E[\delta_\infty(\xi(\eps),\Y_\eps \cup \Y'_\eps) - \delinf|\xi(\eps),\Y_\eps] 
\}
\nonumber \\
\leq 
\E[ (\delta_\infty(\xi(\eps),\Y_\eps \cup \Y'_\eps) + \delinf )^2]^{1/2}
\E[ (\delta_\infty(\xi(\eps),\Y_\eps \cup \Y'_\eps) - \delinf )^2]^{1/2}.
~~~~~
\label{0911}
\eea
By the moments condition (\ref{pbdmom}), the
stabilization condition (\ref{0513a}), and Fatou's lemma, 
$\E[\delinf^4 ] < \infty$. Also, by definition
 $\delinf $ is almost surely the same as $\delta_\infty
(\0, \Po_\la)$ which has the same distribution as 
$\delta_\infty(\xi(\eps),\Po_\la)$  
by translation-invariance, so that
\bean
\E[ \delinf^4 ] =
 \E[ \delta_\infty(\xi(\eps), \Po_\la )^4] 
\\
\geq  
 P[ \Po_\la \cap Q_\0^\eps =\emptyset] 
 \E[ \delta_\infty(\xi(\eps), \Po_\la )^4| \Po_\la \cap Q_\0^\eps =\emptyset] 
\\
=
e^{- \lambda \eps^d} 
 \E[ \delta_\infty(\xi(\eps), \Y_\eps \cup \Y'_\eps )^4], 
\eean
so that
 $\E[ \delta_\infty(\xi(\eps), \Y_\eps \cup \Y'_\eps )^4]$ remains bounded
as  $\eps \downarrow 0$. Combining all these estimates, we obtain
\bea
\limsup_{\eps \downarrow 0}
\E[ (\delta_\infty(\xi(\eps),\Y_\eps \cup \Y'_\eps) + \delinf )^2] < \infty.
\label{0910c}
\eea

As $\eps \downarrow 0$, it is the case that $P[ \Y_\eps \cup \Y'_\eps \neq \Po_\la]$
tends to zero,  and hence
 $$
\delta_\infty(\xi(\eps), \Y_\eps \cup \Y'_\eps)  - \delta_\infty(\xi(\eps), \Po_\la) \toP 0.
$$
Combined with (\ref{0907}) from Lemma \ref{lem0912}, this implies that 
$\delta_\infty(\xi(\eps),\Y_\eps \cup \Y'_\eps) 
\toP \delinf$,
and hence,  using also the fact that $\delta_\infty(\xi(\eps),\Y_\eps \cup \Y'_\eps)$
has uniformly bounded fourth 
moments, we obtain the limit 
$$
\E[ (\delta_\infty(\xi(\eps),\Y_\eps \cup \Y'_\eps) - \delinf )^2] \to 0
{\rm ~~~ as ~~} \eps \downarrow 0.
$$
Combined with (\ref{0911}) and (\ref{0910c}), this
 shows  that
\bea
\lim_{\eps \downarrow 0} \left(
\E \{ ( \E[\delta_\infty(\xi(\eps),\Y_\eps \cup \Y'_\eps)|\xi(\eps),\Y_\eps])^2
- (\E[ \delinf|\xi(\eps),\Y_\eps])^2 \}
\right) =  0.
\label{0912b}
\eea
If we denote by $V_\eps$ the union of the cubes
$Q_z^\eps, z \prec \0$, then by the definition of $Q_z^\eps$
in Section \ref{secnotation}, we find that $V_\eps \subset V_{\eps'}$
for $0 < \eps' < \eps$, and also $\cup_{\eps>0} V_\eps $ is the
half-space $\{(x_1,\ldots,x_d) \in \R^d: x_1 < 0\}$.
Hence the $\sigma$-field generated by $\Y_\eps$ increases
as $\eps$ decreases, and the smallest $\sigma$-field
with respect to which all $\Y_\eps, \eps >0$ are measurable
is the $\sigma$-field $\F$ generated by the Poisson
configuration in the aforementioned half-space (which is the same
 as $\F$ given in the statement of Theorem \ref{1Hclt}).

By the independence of
 $\xi(\eps)$ from  
 $\delinf$ and $\Y_\eps$,
along with the Martingale Convergence Theorem, as $\eps \downarrow 0$
\bean
\E[ \delinf |\xi(\eps),\Y_\eps ] = 
\E[ \delinf |\Y_\eps ] \to \E[\delinf|\F], {\rm ~~~~~ a.s.}
\eean
Since $\E[\delinf^4]< \infty$,
 the variables
$(\E[ \delinf |\xi(\eps),\Y_\eps ])^2 $ are uniformly integrable,
so that
$$
\E[(\E[ \delinf |\xi(\eps), \Y_\eps ])^2] \to \E[(\E[\delinf|\F])^2]
 {\rm ~~as~} \eps \downarrow 0.
$$
Combining this with (\ref{0912b}) and
(\ref{0912c}), we obtain (\ref{0912d})
as required.
$\qed$ \\

To de-Poissonize the limits
(\ref{0516a}), (\ref{0516b}) and obtain (\ref{varlim})  and
(\ref{2clteq}), we use a coupling technique related to that used in
\cite{KL} and \cite{LeeI}. Let $B_0 \in \RR(\R^d)$ with
$|B_0|>0$,  and let $A \in \RR(B_0)$.
Let $U_{1,t},U_{2,t},U_{3,t}\ldots$ be independent and uniformly
distributed over $t B_0$; we assume that
the point processes $\U_{1,t}, \U_{2,t}, \U_{3,t}$ are
coupled by setting
\bea
\U_{m,t} = \{U_{1,t},\dots,U_{m,t}\} , ~~~m \in \N.
\label{couplunif}
\eea
 With this coupling, given point process
set functions $h,h'$, we make the  definition
\bea
R_{m,t}(A) =
\HH(\U_{m+1,t},tA) - \HH( \U_{m,t},tA), \label{Rdef}
\\
R'_{m,t}(A) =
\HH'(\U_{m+1,t},tA) - \HH'( \U_{m,t},tA). \label{Rdef2}
\eea
Let $\la >0$, and
recall from (\ref{mudef}) the definition  
$$
\mu_{\la,t}:= \lambda t^d|B_0|.
$$
We shall use the following coupling lemma, which resembles Lemma 4.2
of \cite{PY}.
\begin{lemm}
\label{couplem}
Suppose $h$  is a point process set function which is 
 strongly stabilizing at intensity $\la$ (i.e., satisfies
(\ref{sstab1}) and (\ref{sstab2}))  with
stabilizing limit $\delinf$.
Suppose $h'$  is a point process set function which is 
 also strongly stabilizing at intensity $\la$
with stabilizing limit $\delinfd$.
Let the
random $d$-vector $Y$ be uniformly distributed over $B_0$
 and independent of $\Po_\la$.
Let $\eps >0$. Then there exists $\eta >0$ and
 $t_0 \geq 1$ such that for all $t \geq t_0$ and all integer $m,m' \in
[(1-\eta)\mu_t, (1+\eta)\mu_t]$
 with $m< m'$, there exists a coupled family of
variables
$D,D',R,R'$
with following properties:
\begin{itemize}
\item
$D$  has the same distribution as $\delinf {\bf 1}_{\{Y \in A\}}$;
\item
$D'$ has the same distribution as $\delinfd {\bf 1}_{\{Y \in A\}}$;
\item
$D$ and $D'$ are independent;
\item
$(R,R')$ have  the same joint  distribution as $(R_{m,t}(A), R'_{m',t}(A))$;
\item
$
P[  \{D \neq R \} \cup \{ D' \neq R' \} ] < \eps.
$
\end{itemize}

\end{lemm}
{\em Proof.} Suppose we are given $t$.  On a suitable probability
space, let
 $\Po$ and $\Po'$ be independent homogeneous Poisson point processes in
 $\R^d$ of intensity $\lambda$; let $U,U',V_1,V_2,\ldots$ be independent
variables uniformly distributed over $tB_0$, independent of $\Po$ and
$\Po'$. 

Let $\Po''$ be the point process consisting of those points
of $\Po$ which lie closer to $U$ than to $U'$ (in the Euclidean norm),
together with
 those points
of $\Po'$ which lie closer to $U'$ than to $U$.  Then $\Po'' $
is a homogeneous Poisson process
 of intensity $\lambda$ on $\R^d$, and moreover it
is independent of $U$ and of $U'$.

Let $N$ denote the number of points of $\Po''$ lying in $t B_0$ (a
Poisson variable with mean
 $\mu_t$).
 Choose an ordering on the points
of $\Po''$ lying in $tB_0$, uniformly at random from all $N!$
possible such orderings.
 Use this ordering to list the points of
$\Po''$ in $tB_0$ as $W_1,W_2,\ldots,W_N$. Also, set $W_{N+1} =
V_1, W_{N+2} = V_2,W_{N+3} = V_3$ and so on.
Define the point process 
$
\W_n := \{W_1, \ldots, W_n\}
$
(for each $n \geq 1$), and the increments
$$
R := \HH(\W_m \cup \{U\},tA) - \HH(\W_m,tA); 
$$
$$
R' = \HH'(\W_{m'-1}\cup \{U,U'\},tA) - \HH'(\W_{m'-1}\cup \{U\},tA) .
$$
The variables  $U,U',W_1,W_2,W_3,\ldots,$ are
independent uniformly distributed variables on $tB_0$, and therefore
the pairs
$(R,R')$ and $(R_{m,t}(A),R'_{m',t}(A))$ have
 the same joint distribution
as claimed.

Let $\tilde{\Po}$ be the translated point process
$\tau_{-U}(\Po )$.
 Similarly, let
 $\tilde{\Po}' : = \tau_{-U'}(\Po')$.
 Then $\tilde{\Po}$ and
 $\tilde{\Po}'$
are independent homogeneous Poisson processes of  intensity $\lambda$ on $\R^d$.
Moreover, $U$ and $U'$ are independent of
 $\tilde{\Po}$
 and $\tilde{\Po}'$.
Let  $S$ be   a  radius of stabilization of $h$ with respect to
 $\tilde{\Po}$, and
let  $S'$ be   a radius of stabilization of $h'$ with respect to
 $\tilde{\Po}'$.
Recall the definition of the add one cost 
$\delta(A,\X)$ at (\ref{add1c}), and
define
$$
D = \delta(B_S(\0), \tilde{\Po} \cap B_S(\0) ) 
{\bf 1}_{\{U \in t A\}}; ~~~~~~~~
D' = \delta'(B_{S'}(\0), \tilde{\Po}'
 \cap B_{S'}(\0)  ){\bf 1}_{\{U' \in t A\}}.
$$
Then $D$ and $D'$ are independent, 
and $D$ has the same distribution as
 $\delinf {\bf 1}_{\{Y \in A\}}$, while
 $D'$ has the same distribution as
 $\delinfd {\bf 1}_{\{Y \in A\}}$.

It remains to show that $(D,D') = (R,R')$ with high probability.
Choose $K$ such that $P[S>K ] < \eps/9$
 and $P[S'>K ] < \eps/9$.
Using the assumption that $B_0$ and $A$ are Riemann measurable,
take $t$ to be
so large that except on an event (denoted $E_0$) of probability
less than $\eps/9$, the positions of $U $ and $U'$ are
Euclidean distance at least $2K$ from $\partial(tB_0)$, from
$\partial(tA)$,
and from each other.
Set $\eta = \eps (2K)^{-d} /(18 \lambda)$.
 We assume $|m- \mu_t| \leq \eta \mu_t$
and $|m'- \mu_t|\leq \eta \mu_t$. Define events $E_1, E, E'$ by 
$$
E_1:= \{|N-m| > 2 \eta \mu_t \}
\cup 
\{|N-m'| > 2 \eta \mu_t \},
$$
$$
E : = \{ \W_m \cap B_K(U) \neq \Po'' \cap tB_0 \cap B_K(U) \}, 
$$
$$
E' : = \{ \W_{m'} \cap B_K(U') \neq \Po''  \cap tB_0 \cap B_K(U') \}.
$$
Event $E$  occurs either if one or more of the
$(N-m)^+$ ``discarded'' points of $\Po''$ 
 lies in $B_K(U)$,
or if one or more of the
$(m-N)^+$ ``added'' points of $\{V_1,V_2,\ldots\}$
 lies in $B_K(U)$,
 and similarly for $E'$.  Hence,
$$
P[E|E_1^{\rm c}] \leq (2 \eta \mu_t) (2K)^d (\lambda/\mu_t) < \eps/9;
~~~
P[E'|E_1^{\rm c}] <  \eps/9.
$$

Using the defining properties 
(\ref{sstab1}) and (\ref{sstab2}) of
the radii of (strong) stabilization $S,S'$
for $\tPo$ and $\tPo'$, and using Boole's inequality,
we obtain for large enough $t$ that
\bean
P[(D,D') \neq (R,R') ] \leq
P[E_0] + P[E_1] + P[S >  K] + P[S' > K] \\
 + P[E \setminus E_1] + P[ E' \setminus E_1]
< \eps. ~~~~~~~~~~~~\qed
\eean

\begin{lemm}
\label{fordepolem}
Let $\la >0$ and let $B_0 \in \RR(\R^d)$ with $|B_0|>0$.
Suppose that $\HH$ and $\HH'$ are point process set functions
which are strongly stabilizing at intensity $\la$,
with stabilizing limit $\delinf$, $\delinfd$ respectively,
and
 $\HH$ and $\HH'$ both satisfy the moments condition $(\ref{ubdmom})$.
 Suppose $A \in \RR(B_0)$. 
Let $g: [1,\infty) \to (0,\infty)$ be a
 function  with $g(t)/t^d \to 0$ as $t \to \infty$.
Then $R_{m,t}$ and $R'_{m',t}$ defined at (\ref{Rdef}), (\ref{Rdef2}) satisfy
\begin{equation}
\lim_{t \to \infty}  \sup_{\mu_{\la,t} -g(t)  \leq m  \leq \mu_{\la,t}+ g(t) }
\left|\E R_{m,t}(A) -
\left( \frac{|A|}{|B_0|} \right)
 \E \delinf
 \right| =0.
\label{depoE}
\end{equation}
Also,
\begin{equation}
\lim_{t \to \infty}  \sup_{\mu_{\la,t} -g(t)  \leq m < m'  \leq \mu_{\la,t} + g(t) }
\left|\E R_{m,t}(A) R'_{m',t}(A) -  
\left( \frac{|A|}{|B_0|} \right)^2 (\E \delinf)\E \delinfd     
\right| =0,
\label{depoE12}
\end{equation}
and
\begin{equation}
\lim_{t \to \infty}  \sup_{\mu_{\la,t} -g(t)  \leq m < m'   \leq \mu_{\la,t} + g(t) }
( \max(\E[ R_{m,t}(A)^4], \E[ R'_{m',t}(A)^4] )   < \infty.
\label{depoE2}
\end{equation}
\end{lemm}
{\em Proof.}
 We start with (\ref{depoE2}); this follows from
the moments condition (\ref{ubdmom}).

Suppose $(t(n), n \geq 1)$ is an arbitrary $(0,\infty)$-valued
sequence tending to infinity as $n \to \infty$. Suppose 
$(m(n), n \geq 1)$
and
$(m'(n), n \geq 1)$ are
 $\N$-valued sequences which satisfy 
\bea
\mu_{\la,t(n)}  -g(t(n)) \leq m(n)
< m'(n) \leq \mu_{\la,t(n)}  + g(t(n)).
\label{0909b}
\eea
By Lemma \ref{couplem}, with $(D,D')$ distributed as in that result
we have as $n \to \infty$ that 
\bea
R_{m(n),t(n)}(A) \tod D; ~~~ 
R_{m(n),t(n)} (A) R'_{m'(n),t(n)} (A) \tod DD'.
\label{0909a}
\eea
By (\ref{depoE2}), the random variables 
$ R_{m(n),t(n)}(A), n \geq 1$, are uniformly integrable, and so are
the variables 
$R_{m(n),t(n)} (A) R'_{m'(n),t(n)} (A) , n \geq 1$.
Hence we have the convergence of expectations corresponding
to the convergence in distribution given by (\ref{0909a}), i.e.,
as $n \to \infty$ we have that
\bean
\E[ R_{m(n),t(n)}(A) ] \to \E[D] =  \frac{|A|}{|B_0|} \E \delinf;  \\ 
\E [ R_{m(n),t(n)} (A) R'_{m'(n),t(n)} (A) ]
 \to   \E[DD'] =\frac{|A|^2}{|B_0|^2} (\E \delinf ) \E \delinfd  ,
\eean
and since the choice of $t(n)$, $m(n)$ and  $m'(n)$ was arbitrary 
subject to $\limn (t(n)) = \infty$ and
to (\ref{0909b}),  this gives us (\ref{depoE}) and (\ref{depoE12}).
$\qed$

\vskip 1em
\noindent
{\em Proof of Theorem \ref{2Hclt}.} 
Assume $(t_n)_{n \geq 1}$ is a $(1,\infty)$-valued sequence
satisfying (\ref{tncond}), which says that
$\la t_n^d |B_0| -n $ is $O(n^{1/2})$ as $n \to \infty$.

Assume the point processes 
$\Po_{t_n},$ 
$ \U_{1,t_n}, $ $ \U_{2,t_n},$ $\U_{3,t_n}, \ldots$
are  coupled  by having $\U_{m,t_n}$ defined by (\ref{couplunif})
and
setting $\Po_{t_n} = \{U_{1,t_n}, U_{2,t_n},\ldots,U_{N_n,t_n} \}$ with
$N_n$ an independent Poisson variable with mean
 $\mu_n := \mu_{\la,t_n} = \lambda t_n^d |B_0|$.
For $1 \leq j \leq k$, let
\bean
\HHH_n^j : =  \HH^j(\U_{n,t_n},t_n A_j) ; ~~~~
\tilde{\HHH}_n^j :=   \HH^j(\Po_{t_n},t_nA_j).
\eean
Define the $k$-vector 
$$
\alpha  : = (\alpha_j)_{j=1}^k, ~~~~{\rm with~~} \alpha_j :=
 \left( \frac{|A_j|}{|B_0|}\right)\E[\delinfj].
$$ 
The first step is to prove that
as $n \to \infty$,
\begin{equation}
\E \left[ (n^{-1/2} (\tilde{\HHH}^j_n - \HHH^j_n -
 (N_n-n) \alpha_j))^2 \right]
 \to 0.
\label{fordepo}
\end{equation}
To prove this, (writing $t(n)$ for $t_n$ when typographically convenient),
 note that the expectation in the left hand  side
is equal to
\bear
\sum_{m: |m-\mu_{n} |\leq n^{3/4} } 
 \E \left[ n^{-1}
 \left(
\HH^j(\U_{m,t(n)},t_n A_j ) - \HH^j(\U_{n,t(n)},t_nA_j) - (m-n)
  \alpha_j \right)^2 \right] P[N_n=m]
\nonumber
\\
 + n^{-1} \E \left[  \left(
 \tilde{\HHH}^j_n - \HHH^j_n  - (N_n- n) \alpha_j \right)^2
 {\bf 1}\{ |N_n - \mu_{t(n)}|
> n^{3/4} \} \right]. ~~~~~~~~~~~~~~~
\label{fordepo1}
\eear
Let $\eps >0$.
By (\ref{Rdef}) and
Lemma \ref{fordepolem}, there exists $c>0$ such that for
large enough $n$
and all $m$ with $n \leq m \leq \mu_{n}+ n^{3/4}$,
\bean
\E[ (\HH^j(\U_{m,t(n)},t_nA_j) - \HH^j (\U_{n,t(n)},t_nA_j ) -
 (m-n) \alpha_j )^2] ~~~~~~
\\
= \E \left[  
\left( \sum_{\ell =n}^{m-1}
(R^j_{\ell ,t(n)}(A_j)- \alpha_j) \right)^2 \right]
\leq \eps(m-n)^2 + c(m-n),
\eean
where the bound comes from expanding out the double sum
 arising from the expectation of the squared sum;
the $c(m-n)$ term comes from bounding the diagonal terms
using (\ref{depoE2}) and the fact that bounded fourth moments imply
bounded second moments.
A similar argument applies when $\mu_{n}-n^{3/4} \leq m \leq n$, and
hence
the first term in (\ref{fordepo1}) is bounded by the expression
\bean
n^{-1} \E[ \eps(N_n-n)^2 + c |N_n-n| ] = n^{-1}  [\eps ( 
\E[(N_n-\mu_n)^2 ] + (\mu_n-n)^2)  + c  \E[ |N_n-n| ]  ] . 
\eean
By assumption (\ref{tncond}), we have that $\mu_n \sim n$ and
 $c' : = \limsup (\mu_n-n)^2/n < \infty$, so that
for large $n$  the first term in (\ref{fordepo1}) is bounded by
 $ ( 3 + c' )\eps$ 
 for $n$ large enough. By the
uniform bound (\ref{polybd})
 and the Cauchy-Schwarz inequality,
there  is a  constant $\beta_3$  such that the
 second term in (\ref{fordepo1}) is bounded by
$ \beta_3 n^{\beta_3} (P [|N_n-\mu_{n}| > n^{3/4}])^{1/2} $,
which tends to zero, e.g. by Lemma 1.4 of \cite{Pbk}.
 Since $\eps$ is arbitrary and does not depend on $c'$,
this completes the proof of (\ref{fordepo}).

Let  $b_1, \ldots, b_k$ be arbitrary real constants. 
Define the column vector
$\bb := (b_1,\ldots, b_k)'$.
Let
\bean
\HHH_n : = \sum_{j=1}^k b_j \HHH_n^j; ~~~~
\HHH'_n = \sum_{j=1}^k b_j \tilde{\HHH}_n^j.
\eean
We prove convergence of $n^{-1}\Var(\HHH_n)$, using the identity
$$
n^{-1/2} \HHH'_n = n^{-1/2} \HHH_n + n^{-1/2} (N_n-n)\balpha'\bb
+ n^{-1/2} (\HHH'_n- \HHH_n - (N_n-n)\balpha'\bb).
$$
In the right hand side, the third term has variance tending to zero
by (\ref{fordepo}), while the second term has variance tending to
 $(\balpha'\bb)^2$ and is independent of the first term. It follows that
with the matrix $\Sigma^\la =(\sigma_{ij}^\la)_{i,j=1}^k$
 given by Theorem \ref{1Hclt},
and the matrix $\Sigma^{\la,A} = (\sigma_{ij}^{\la,A})_{i,j=1}^k$ given by
$$
\sigma_{ij}^{\la,A} := \frac{\sigma_{ij}^\la|A_i \cap A_j| }{|B_0|},
$$
we have from Theorem \ref{1Hclt} that
$$
\bb' \Sigma^{\la,A} \bb =
\limn n^{-1} \Var(\HHH'_n) = \limn ( n^{-1} \Var(\HHH_n) )+  (\alpha'  \bb)^2,
$$
so that $\Sigma^{\la,A} - \balpha \balpha'$ is nonnegative definite and
 $n^{-1} \Var(\HHH_n) \to
\bb' (\Sigma^{\la,A} - \balpha \balpha')\bb
 $.
This gives us (\ref{varlim}).

The proof of Theorem \ref{1Hclt} (since it 
is derived by taking linear combinations) tells us that
$
n^{-1/2}
(\HHH'_n - \E \HHH'_n )\tod
\NN(0,\bb' \Sigma^{\la,A} \bb).
$
Combined with (\ref{fordepo}) this gives us
\bea
n^{-1/2}(\HHH_n - \E \HHH'_n + (N_n-n) \alpha' \bb) \tod
\NN(0,\bb' \Sigma^{\la,A} \bb).
\label{1123}
\eea
Recall that $\mu_n := \lambda t_n^d |B_0| = \E N_n$.
Since $n^{-1/2}(N_n-\mu_n)\balpha' \bb  $ is independent of
$\HHH_n$ and is asymptotically normal with
mean zero and variance $(\balpha' \bb)^2 = \bb' \alpha \alpha'\bb $,
we can deduce from (\ref{1123}), by considering
characteristic functions, that
\begin{equation}
n^{-1/2}(\HHH_n - \E \HHH'_n + (\mu_n-n)\balpha'\bb) \tod
\NN(0,\bb' (\Sigma^A -\balpha \balpha')\bb).
\label{fordepo2}
\end{equation}
By (\ref{fordepo}), the expectation of
$n^{-1/2} (\HHH'_n - \HHH_n  - (N_n-n) \balpha' \bb)$ tends to zero,
 so in (\ref{fordepo2})
we can replace $-\E \HHH'_n+ (\mu_n -n)\balpha'\bb$ by $-\E \HHH_n$, 
which gives us
$$
n^{-1/2} (\HHH_n - \E \HHH_n) \tod \NN 
(0, \bb'(\Sigma^A - \balpha \balpha') \bb  ).  
$$
Then (\ref{2clteq}) follows by the Cram\'er-Wold device. $\qed$ \\

\begin{lemm}
Let $\la >0$.
Suppose the graph $G:= G(\X)$, defined
for each locally finite $\X \subset \R^d$,
 is translation invariant and satisfies
the stabilization conditions (\ref{gstab1}) and (\ref{gstab2}).
Then with probability 1, 
for each $X \in \Po_\la^0$ there exists 
$R(X)<\infty$ such that the set of edges of $G(\Po_\la^0)$
incident to $X$ is unaffected by changes to $\Po_\la^0$
outside $B_{R(X)}(X)$,   
for each $X \in \Po_\la$ there exists 
$R(X)<\infty$ such that the set of edges of $G(\Po_\la)$
incident to $X$ is unaffected by changes to $\Po_\la$
outside $B_{R(X)}(X)$.
\label{PY4lem}
\end{lemm}
{\em Proof.} The existence of finite $R(X)$ for all
$X \in \Po_\la^0$  is given by Lemma 3.3
 of \cite{PY4}. The existence of finite $R'(X)$
for all $X \in \Po_\la$  is proved in the course of the
proof of Lemma 3.3 of \cite{PY4}. $\qed$ \\

\noindent
{\em Proof of Lemma \ref{lem0516}.}  
First set $h(\X,A)= L^G_\phi(\X,A)$, as defined at (\ref{phiw1}).
Let the random variable $R$ satisfy (\ref{gstab1}) and (\ref{gstab2}). 
Let $A \in \RR(\R^d)$, and let
 ${\cal A} \subset \R^d \setminus B_R(\0)$ be finite.
Then if $B_R(\0) \subseteq A$,  the increment
$\delta(A, (\Po_\la \cap B_R(\0)) \cup {\cal A})$
is equal to
\bea
 \left(
 \sum_{e \in \Ed^+(\0;\Po_\la^0 \cap B_R(\0))} \phi(|e|)
\right)
- 
\sum_{e \in \Ed^-(\0;\Po_\la^0 \cap B_R(\0))} \phi(|e|) , 
\label{0513c}
\eea
since all added and removed edges  have both endpoints in $A$. 
Hence (\ref{sstab1}) holds  with $\delinf $ equal
to the expression displayed in (\ref{0513c}). 
If instead $A \cap B_R(\0) = \emptyset$, then 
$\delta(A, (\Po_\la \cap B_R(\0) ) \cup {\cal A} ) =0$ 
since added and removed edges
have neither endpoint in $A$. Hence, (\ref{sstab2}) holds.

Next, suppose we set $h(\X,A)= V^G_\psi(\X,A)$,
where $\psi \in B(\K_\ka)$, with $\ka \in \N$.
We assert that there exists an almost surely finite
random variable $R$ such that
(\ref{gstab1}) and (\ref{gstab2}) hold,
and such that 
for every vertex $X$ of $G(\Po_\la^0 \cap B_{R}(\0))$ at a graph distance
at most $2 \kappa $ from 
 some endpoint of some
edge in either
 ${\cal E}^+(\0; \Po_\la^0 \cap B_{R}(\0))$
or
 ${\cal E}^-(\0; \Po_\la^0 \cap B_{R}(\0))$,
the set of edges incident to $X$
is unaffected by changes outside $B_{R}(\0)$.
The existence of such an $R$  follows from Lemma \ref{PY4lem} along
with an inductive argument in $\kappa$.

Let ${\cal A} \subset \R^d \setminus B_{R}(\0) $ 
be finite.  Suppose $X \in (\Po_\la \cap B_R(\0) ) \cup {\cal A}) $
lies at a graph distance
more than $\kappa$ 
 in $G((\Po_\la^0 \cap B_{R}(\0)) \cup {\cal A})$ 
from
 any endpoint of any
edge in
either 
 ${\cal E}^+(\0; (\Po_\la^0 \cap B_R(\0)) \cup{\cal A})$ 
or
 ${\cal E}^-(\0; (\Po_\la^0 \cap B_R(\0)) \cup{\cal A})$ 
(all vertices in ${\cal A}$ fall in this category). 
Then
$$
\psi(G_{X,\ka}[(\Po_\la^0 \cap B_R(\0) ) \cup {\cal A}]) =
 \psi(G_{X,\ka}[(\Po_\la \cap B_R(\0) ) \cup {\cal A}]).
$$
Also, for the remaining $X \in \Po_\la \cap B_R(\0)$,
 at a graph distance at most $\kappa $ from 
the endpoint some edge in either
 ${\cal E}^+(\0; (\Po_\la^0 \cap B_R(\0) )\cup {\cal A})$
or
 ${\cal E}^-(\0; (\Po_\la^0 \cap B_R(\0) )\cup {\cal A})$,
the value of $\psi(G_{X,\ka}((\Po^0  \cap B_R(\0) )\cup {\cal A})) 
- \psi(G_{X,\ka}((\Po \cap B_R(\0) )\cup {\cal A}))$ 
is unaffected by changes to the set ${\cal A}$  outside $B_{R}(\0)$.
The conditions
 (\ref{sstab1}) and (\ref{sstab2}) (with $S=R$) follow for  this case.

Suppose now that uniqueness of the infinite component holds,
and set $h(\X,A) = K^G(\X,A)$. The stabilization conditions
are proved, essentially by a slight modification of the proof
of Proposition 6.1 of \cite{PY}. For the convenience of the 
reader, we describe the argument in the present, more general
context.

Let $R$ be a radius of stabilization, as given at (\ref{gstab1})
and (\ref{gstab2}).
Choose a finite $R' >R$ such that for any two points of 
$\Po$ in $B_{R}(\0)$,
either they can be connected by a path in $G(\Po)$ all of whose nodes
lie in $B_{R'}(\0)$, or at least one of them lies in a finite 
component contained in  $B_{R'}(\0)$, and such that a similar statement
holds for $\Po^0$. The proof that we can choose such an $R'$
 is based on
the uniqueness of the infinite component in $G(\Po)$ 
and $G(\Po^0)$, and is given in more detail in \cite{PY}.

By Lemma 3.3 of \cite{PY4}, there almost surely exists 
$R''> R'$ such that for all $X \in \Po \cap B_{R'}(\0)$, the
set of edges incident to $X $ in $G(\Po^0)$ is 
unaffected by additions or deletions of 
points outside $B_{R''}$, and moreover, by the proof
of Lemma 3.3 of \cite{PY4}, we can choose $R''$ to 
be so large  that  in addition,
the
set of edges incident to $X $ in $G(\Po)$ is 
unaffected by additions or deletions of 
points outside $B_{R''}$.

Suppose that  $A \in \RR(\R^d)$ and
$B_{R''}(\0) \subseteq  A $. Suppose ${\cal A}$ is disjoint from
$B_{R''}(\0)$.
When we change from
 $G(\Po \cap B_{R''}(\0) \cup {\cal A})$ to
 $G(\Po^0 \cap B_{R''}(\0) \cup {\cal A})$, the effect is first to
add a vertex at the origin, then to add the edges of
$\Ed^+(\0;\Po^0)$, and then to remove the edges of $\Ed^-(\0;\Po^0)$.
Consider adding successive edges, in some specified
order.
 Each edge reduces the number of components
that intersect $A$
by 1 if it joins
two points that were previously not connected by a path, and otherwise
does not affect the number of components. The question
of whether a particular added edge changes the number of components
is determined by the graph structure of 
the restriction of $G(\Po)$   to vertices in $B_{R'}(\0)$, and therefore
does not depend on $A$ or ${\cal A}$ (always presuming $A  \in \RR$ and
 $ B_{R''}(\0) \subseteq A$).
 A similar argument applies with deleted edges.

It follows from the above that if we set 
$ \delinf : = \delta(B_{R''}(\0),\Po \cap B_{R''}(\0))$
and $S = R''$,
then (\ref{sstab1})  holds.

Now suppose that $A \cap B_{R''} (\0) = \emptyset$ (and ${\cal A} $ is
also disjoint from $B_{R''}(\0)$ as before).
Consider again the process of
 successive additions and deletions described above.
If an added edge connects two previously disconnected components,
then  at least one of them has a vertex set entirely contained
in $B_{R'}(\0)$, and therefore does not have any vertices in $A$,
and so this change does not cause any increment in the number
of components that have at least one vertex in $A$. A similar
argument applies with removed vertices;  hence,  
if  $B_{R''} (\0)\cap  A = \emptyset $ we have
  $\delta(A,(\Po \cap B_{R''}(\0)) \cup {\cal A})  =0$, 
so that (\ref{sstab2})  holds.
 $\qed$ \\

\noindent
{\bf Acknowledgement.} I thank the referee for carefully
reading the first version of this paper, and pointing out
some inaccuracies and obscurities therein.

Department of Mathematical Sciences

University of Bath  

Bath  BA2 7AY

United Kingdom

\texttt{M.D.Penrose@bath.ac.uk}

%
%
%

\end{document}